\documentclass[a4paper,11pt]{article}

\usepackage[latin1]{inputenc}
\usepackage{epsfig,color}
\usepackage{amsmath,amsfonts,amssymb}

\usepackage{mathptmx}

\input{macro}

\newlength{\hauteur}
\newcommand{\rond}[1]{\settoheight{\hauteur}{$#1$}\raisebox{\hauteur}{$\scriptscriptstyle\,\,\circ$}\hspace{-1ex}{#1}}%

\renewcommand{\A}{\mathbb{A}}%
\renewcommand{\P}{\mathfrak{P}}%
\newcommand{\intervalle}[1]{\pmb{\lfloor}#1\pmb{\rceil}}
\newcommand{\zhilin}{Zhilinski\'{\i}}
\newcommand{\sadov}{Sadovski\'{\i}}

\title{Moment polytopes\\ for symplectic manifolds with monodromy\\
}

\author{V\~u Ng\d oc San}
\date{}
\begin{document}
\maketitle

\begin{abstract}
  A natural way of generalising Hamiltonian toric manifolds is to
  permit the presence of generic isolated singularities for the moment
  map.  For a class of such ``almost-toric 4-manifolds'' which admits
  a Hamiltonian $S^1$-action we show that one can associate a group of
  convex polygons that generalise the celebrated moment polytopes of
  Atiyah, Guillemin-Sternberg. As an application, we derive a
  Duistermaat-Heckman formula demonstrating a strong effect of the
  possible monodromy of the underlying integrable system.
\end{abstract}

\vfill

\noindent\textbf{Keywords~: } moment polytope, circle action,
semi-toric, Duistermaat-Heckman, monodromy, symplectic geometry,
Lagrangian fibration, completely integrable systems.

\noindent\textbf{Math. Class.~:}
53D05, 
53D20, 
37J15, 
37J35, 
57R45  

\newpage

\section{Introduction}

Let $M$ be a compact connected symplectic manifold, equipped with an
effective Hamiltonian action of a torus $\T^k$. A moment map for this
action is a map $\Phi:M\fleche\RM^k$ (where $\RM^k$ is viewed as the
dual of the Lie algebra of $\T^k$) whose components generate commuting
Hamiltonian flows which are independent almost everywhere and thus
define the given effective $\T^k$ action.  In 1982,
Atiyah~\cite{atiyah-convex} and
Guillemin-Sternberg~\cite{guillemin-sternberg} discovered
independently that the image of $\Phi$ is very special: it is a convex
polytope. This polytope encodes many pieces of information about
$(M,\Phi)$; if the action is completely integrable in the sense that
$2k$ is the dimension of $M$ then Delzant~\cite{delzant} actually
proved that the moment polytope completely determines $(M,\Phi)$,
thereby showing that $M$ is in fact a toric variety.

The theory of Hamiltonian actions on symplectic manifolds has more
recently been extended to include non-compact manifolds, provided the
momentum map is proper. Then all the results essentially persist.

From the point of view of classical mechanics and applications to
quantum mechanics, one is generally more interested in the particular
Hamiltonian function under study than in the underlying manifold.
Toric manifolds are perfectly good phase spaces for many relevant
examples, but the class of toric Hamiltonians or toric momentum maps
is by far too narrow.

Mechanical systems usually will show up more complicated singularities
that those allowed by toric momentum maps. A much more flexible notion
to use instead of completely integrable toric actions is completely
integrable systems, which means that one is given a ``momentum map''
$\Phi=(f_1,\dots,f_n)$ with the only requirement that $\{f_i,f_j\}=0$
for all $i,j$ and $df_1,\dots,df_n$ are independent almost everywhere.
In other words $\Phi$ is a momentum map for a local Hamiltonian action
of $\RM^n$, which is locally free almost everywhere. In this
generality, the image of the momentum map (sometimes called the
bifurcation diagram) is still of great interest but has a much more
complicated structure. Even with the requirement that all
singularities be non-degenerate à la Morse-Bott, the global picture is
much richer than a convex polytope (see for
instance~\cite{fomenko-msri} for 2 degrees of freedom).  Nevertheless,
under the assumption that the momentum map is proper (and a submersion
almost everywhere), the Liouville-Arnold-Mineur theorem (or
action-angle theorem) still says that each regular orbit of $\Phi$ is
an $n$-torus in a neighbourhood of which the action is toric.  Hence
the main question is how to globalise this Liouville-Arnold-Mineur
theorem and has two related facets. First is the study of the
topological invariants of the restriction of the momentum map to
regular points: this was explained in Duistermaat's
paper~\cite{duistermaat}. Secondly one has to study the singularities
of $\Phi$ and how they show up in topological or symplectic
invariants. A global picture for this was developed by Nguyên Tiên
Zung~\cite{zung-II}.

In our paper we bring both theories (toric actions and integrable
systems) together in the sense that we construct moment polytopes with
some of the usual properties (rationality, convexity) for momentum
maps that are \emph{not} toric. Our initial motivation was that these
polytopes happen to be excellent tools for the semiclassical study of
the eigenvalues of quantised Hamiltonians \cite{san-redistribution}.

We deal here with symplectic 4-manifolds endowed with a completely
integrable system $\Phi=(J,H)$, $\{J,H\}=0$, such that $J$ alone is a
proper momentum map for an $S¹$-action on $M$. Such a $\Phi$ will be
called \emph{semi-toric}. Then we will assume that all singularities
are non-degenerate (in the sense of Eliasson) \emph{without hyperbolic
  component}. In other words we allow -- in addition to tori --
singular fibres of focus-focus type, which are pinched tori. A torus
pinched once is an immersion of a sphere with one double point and is
considered as the ``simplest'' singular fibre for a 2-torus fibration
(see~\cite{matsumoto}). We prove the following result

\begin{theo*}[proposition~\ref{prop:fibres} and
  theorem~\ref{theo:image}]~\\
  The image of $\Phi$ is simply connected, $\Phi$ has connected
  fibres, and the critical values of $\Phi$ are exactly the points in
  the boundary of the image, plus a finite number of isolated points
  corresponding to the focus-focus fibres.
\end{theo*}

Then we show that, in spite of the fact that focus-focus fibres imply
non-trivial monodromy and hence the impossibility of constructing
global action variables and a $\T^2$-action, one can naturally
transform the image of $\Phi$ into a rational convex polygon which is
almost everywhere the image of a (local) momentum map for a 2-torus
action with the same foliation by tori as $\Phi$. This is the content
of theorem~\ref{theo:polygon}. Such generalised ``moment polytopes''
are not unique; on the contrary the set of all possible polytopes for
a given system has a natural structure of an abelian group isomorphic
to $(\ZM/2\ZM)^{m_f}$, where $m_f$ is the number of focus-focus fibres
(proposition~\ref{prop:group}).

This construction will be used to give a simple formula for the
Duistermaat-Heckman function associated to the $S^1$ action generated
by $J$, which shows clearly the role played by the possible
\emph{monodromy} of the integrable system.

\begin{theo*}[theorem~\ref{theo:DH}]
  If $\alpha^+(x)$ (resp. $\alpha^-(x)$) denotes the slope of the top
  (resp. bottom) boundary of a generalised moment polytope for $\Phi$,
  then the derivative of the Duistermaat-Heckman function is
  \[
  \rho'_J(x)=\alpha^+(x)-\alpha^-(x)
  \]
  and is piecewise constant on $J(M)$. Discontinuities appear at the
  absciss\ae{} $x$ of critical values of $\Phi$ of maximal corank and
  are given by the jump formula~:
  \begin{equation}
    \rho'_J(x+0)-\rho'_J(x-0)= -k(x) - e^+(x) - e^-(x),
    \label{equ:DH-intro}
  \end{equation}
  where $k(x)\in\NM^*$ is an associated monodromy index, and
  $e^{\pm}(x)$ are non-negative contribution of corners of the
  polytope, of the form

  \[
  e^\pm = - \frac{1}{a^\pm b^\pm} \geq 0,
  \]
  where $a^\pm$, $b^\pm$ are the isotropy weights for the $S¹$ action
  at the corresponding vertices.
\end{theo*}

Since quantities in the right-hand side of~\eqref{equ:DH-intro} are
negative, we see that singularities --- and especially those inducing
monodromy --- have a strong effect on the geometry of the polygon. In
particular this yields, as a corollary, the striking result:

\begin{theo*}[corollary~\ref{coro:compact}]
  If $M$ admits a semi-toric momentum map $\Phi=(J,H)$ with at least
  focus-focus critical fibres and such that $J$ has a unique minimum
  (or maximum) then $M$ is compact.
\end{theo*}

\section{Almost toric momentum maps}

Before studying semi-toric momentum map we shall need some general
results about a wider class of integrable systems which are not far
from defining a torus action on $M$, in a suitable sense. The main
result of this section which will be crucial for our purposes is the
description of the image of the moment map for such ``almost-toric''
systems, when the fibres are connected
(proposition~\ref{prop:fibres}).

Although we are interested here in two degrees of freedom, the results
can probably be extended (\emph{mutatis mutandis}) to an arbitrary
dimension.

Let $M$ be a connected symplectic 4-manifold, and $(J,H)$ a completely
integrable system on $M$: $\{J,H\}=0$, such that
$\Phi:=(J,H):M\fleche\RM^2$ is a proper map.

\begin{defi}
  \label{defi:toric}
  A proper $\Phi$ will be said of \textbf{toric type} if there exists
  an effective, completely integrable Hamiltonian $\T^2$-action on $M$
  whose momentum map is of the form $F=f\circ \Phi$, where $f$ is a
  local diffeomorphism on the image of $\Phi$.
\end{defi}
The topology and the symplectic geometry of Hamiltonian $\T²$-actions
are a classical subject, described by what we call the \emph{convexity
  theorem} by Atiyah~\cite{atiyah-convex},
Guillemin-Sternberg~\cite{guillemin-sternberg}, the
\emph{connectedness theorem}, which is generally tied to the
former~\cite{atiyah-convex}, and the \emph{uniqueness theorem} by
Delzant~\cite{delzant}. Note that these results have been generalised
for non-compact manifolds in case of proper momentum maps by Lerman \&
al~\cite{lerman-al}). We will use in this work the following
statements:
\begin{theo}[\cite{lerman-al}]
  If $F$ is a proper momentum map for a Hamiltonian $\T^k$-action on a
  symplectic manifold $M$ then
  \begin{itemize}
  \item the fibres of $F$ are connected;
  \item the image of $F$ is a rational convex polyhedron
  \end{itemize}
\end{theo}
A rational convex polyhedron is by definition a set which can be
obtained near each point by a finite intersection of closed
half-spaces whose boundary hyperplanes admit normal vectors with
integer coefficients.
\begin{prop}
  \label{prop:diffeo}
  In the definition above, $f$ is a diffeomorphism from the image of
  $\Phi$ into the image of $F$. Therefore the fibres of $\Phi$ are
  connected.
\end{prop}
\begin{demo}
  $f$ is surjective by definition. Let us show that it is injective.
  Let $c$ in the image of $F$. Since $F^{-1}(c)$ is connected and
  $F^{-1}(c)=\Phi^{-1}(f^{-1}(c))$, $f^{-1}(c)$ must be connected.
  Since $f$ is a local diffeomorphism, $f^{-1}(c)$ is just a point;
  hence $f$ is injective.
\end{demo}
\begin{rema}
  A weaker definition would be that there exists an effective,
  completely integrable $\T^2$-action on $M$ which leaves $\Phi$
  invariant. This is indeed strictly weaker since this would allow
  $\Phi=g\circ F$ where $g$ is any local diffeomorphism (=immersion),
  but not necessarily a global one (for instance $g$ can send a square
  to an annulus). See also proposition~\ref{prop:simply-connected}
  below.
\end{rema}

We shall be interested here in momentum maps that sometimes fail to be
of toric type.
\begin{defi}
  A proper $\Phi$ is called \textbf{almost-toric} if all the
  singularities are non-degenerate in the sense of Eliasson without
  hyperbolic blocks.
\end{defi}
Note that Symington~\cite{symington-four} has independently introduced
the same definition, and discussed many of its consequences of
topological nature. For a discussion and references on the notion of
Eliasson's non-degeneracy condition, see for instance \cite{san-fn}.
At a critical point of rank zero ($d\Phi(m)=0$) this means that a
generic linear combination of the Hessians $J''(m)$ and $H''(m)$
defines a Hamiltonian matrix (via multiplication by the linearised
symplectic form) that has pairwise distinct eigenvalues. Then
Eliasson's theorem says that the Lagrangian foliation near such a
critical point can be linearised in the $\Cinf$ category.
\begin{prop}[\cite{delzant}]
  \label{prop:elliptic}
  If $\Phi$ is of toric type then $\Phi$ is almost-toric (with only
  elliptic singularities).
\end{prop}
\begin{demo}
  This is a standard argument.  Let $F$ be a momentum map for the
  $\T^2$-action. By definition the singularities of $\Phi$ are the
  same as those of $F$. Now the result follows from the fact that a
  torus action is linearisable near a fixed point. Details can be
  found for instance in~\cite{delzant}
\end{demo}
\begin{prop}
  \label{prop:connected}
  If all the singularities of $\Phi$ are non-degenerate and the set of
  regular values of $\Phi$ is connected then $\Phi$ is almost-toric.
\end{prop}
\begin{demo}
  If a singular point of $\Phi$ has a hyperbolic block, then because
  of the normal form for non-degenerate singularities, there is an
  embedded line segment of critical values in the interior of the
  image of $\Phi$. We conclude by the following lemma.
\end{demo}
\begin{lemm}
  \label{lemm:segment}
  Assume all the singularities of $\Phi$ are non-degenerate. If there
  is an embedded line segment of critical values in the interior of
  the image of $\Phi$, then the set $B_r$ of regular values of $\Phi$
  is not connected.
\end{lemm}
\begin{demo}
  Let $\gamma$ be this segment. Choose an orientation in $\RM^2$ and
  along $\gamma$: since the set $B_r$ of regular values is open and
  dense in $\Phi(M)$, there exists small disjoint open balls on each
  side of $\gamma$. Because all singularities are non-degenerate,
  $\gamma$ can be extended (in both directions) until it reaches a
  singular value of rank zero: elliptic-elliptic, hyperbolic-elliptic
  of hyperbolic-hyperbolic. In all cases $\gamma$ is connected to one
  or several other branches of critical values.  Choose one
  arbitrarily, and continue forever (in both directions).  Since
  $\Phi$ is proper the set of critical values is compact in any
  compact of $\RM^2$, therefore only two things can happen: either
  $\gamma$ intersects itself, or $\gamma$ goes to infinity (goes out
  of any compact) in both directions. In both cases $\gamma$
  disconnects $B_r$.
\end{demo}

In general fibres of almost-toric momentum maps need no be connected.
For instance if $F$ is a toric momentum map and $f$ is a non-injective
immersion of the image of $F$ into $\RM^2$, then $f\circ F$ is
almost-toric with non-connected fibres. However we have the important
proposition below:
\begin{prop}
  \label{prop:fibres}
  Assume $\Phi$ is almost-toric. Consider the following statements:
  \begin{enumerate}
  \item The fibres of $\Phi$ are connected;
  \item the set $B_r$ or regular values of $\Phi$ is connected;
  \item $B_r$ is ``locally connected'': for any value $c$ of $\Phi$,
    for any sufficiently small ball $D$ centred at $c$, $B_r\cap D$
    is connected;
  \item $B_r=\rond{B}\setminus\{c_1,\dots,c_{m_f}\}$, where
    $B=\Phi(M)$, $m_f\leq\infty$ and $c_j$'s are the (isolated) values
    by $\Phi$ of the focus-focus singularities.
  \end{enumerate}
  Then we have $1\impliq 2$, and $2$, $3$, $4$ are equivalent.
\end{prop}
\begin{demo}
  Recall that if $c$ is a critical value of $\Phi$, we call
  $\Phi^{-1}(c)$ a critical fibre. Sometimes we say also a singular
  fibre.
  \paragraph{$1\impliq 2$ :}
  Since $\Phi$ is almost-toric, the singular fibres are either points
  (elliptic-elliptic), circles (codimension 1 elliptic) or pinched
  tori (focus-focus). They do not include regular tori since the
  fibres are assumed to be connected. Only codimension 1 elliptic
  critical values can appear in 1-dimensional families, and
  elliptic-elliptic critical values appear at the end of these
  families. Focus-focus pinched tori are isolated. Therefore the union
  of all critical fibres is a locally finite union of points,
  cylinders and pinched tori, and therefore of codimension 2. Hence
  the complementary set is connected, and therefore its image by
  $\Phi$ also.

  \paragraph{$2\impliq 3$ :}
  Because of the normal forms of the singularities, the only way to
  disconnect a small disc $D\subset \Phi(M)$ is by an embedded segment
  of critical values. But then $B_r$ would not be connected by
  Lemma~\ref{lemm:segment}.

  \paragraph{$3\impliq 4$ :}
  If there is a critical value $c$ in the interior of $B$, then it is
  either isolated (then it must be the image of a focus-focus point)
  or inside an embedded line segment of critical values (which would
  come from codimension 1 elliptic singularities). But the latter case
  is obviously in contradiction with the hypothesis of local
  connectedness.

  \paragraph{$4\impliq 2$ :}
  $B$ is pathwise connected since $M$ is a connected manifold. Suppose
  $c$ and $c'$ are in $\rond{B}$. They can be connected by a path in
  $B$. If this path meets the boundary $\partial B$ (recall that $B$
  is closed since $F$ is proper), it can be pushed inside $\rond{B}$
  using the normal form of elliptic singularities.  Hence $\rond{B}$
  is connected, and the result follows.
\end{demo}
\begin{rema}
  In the proposition above, $2\impliq 1$ is not true. One can imagine
  a torus bundle over an annulus, where the fibre consists of two
  $2$-dimensional tori which swap when going round the annulus.  Note
  however that $2\impliq 1$ is true in case $B_r$ is simply connected,
  as shown in Proposition~~\ref{prop:simply-connected} below. One
  might also conjecture that it is true also when $B$ is simply
  connected.
\end{rema}
\begin{rema}
  The points $c_i$ are called \emph{nodes} in the terminology
  introduced by Symington~\cite{symington-four}.
\end{rema}

In the next section, moment polyhedrons will be defined for some
almost-toric actions. This would happen obviously if the action were
actually toric:
\begin{prop}
  \label{prop:simply-connected}
  If $\Phi$ is almost-toric then $\Phi$ is of toric type if and only
  if the set of regular values of $\Phi$ is connected and simply
  connected.
\end{prop}
\begin{demo}
  Assume $B_r$ is connected and simply connected.  Using the
  connectedness we know from point $4$ of
  proposition~\ref{prop:connected} that
  $B_r=\rond{B}\setminus\{c_1,\dots,c_{m_f}\}$. By the simple
  connectedness we must have $m_f=0$. The fibres corresponding to the
  values in the boundary $\partial B$ can only contain
  elliptic-elliptic fixed points and codimension 1 elliptic circles
  (otherwise $\Phi$ would take values in a small ball centred at our
  point in the boundary...). Therefore the union of all these fibres
  is of codimension 2 so $\Phi^{-1}(B_r)=\Phi^{-1}(\rond{B})$ is
  connected.
  
  Now, since $\pi_0(B_r)=1$ and $\pi_1(B_r)=1$ , the homotopy sequence
  of the fibration $\Phi_{\restr \Phi^{-1}(B_r)}$, implies that
  $\pi_0(\Phi^{-1}(B_r))\simeq\pi_0(\mathcal{F})$, where $\mathcal{F}$
  is the generic fibre of $\Phi$. Hence $\pi_0(\mathcal{F})=1$: the
  fibres are connected.

  Now let $B\subset\RM^2$ be the image of $\Phi$. For each $c\in B$ we
  define the $\ZM$-module of germs of basic action variables at $c$,
  ie germs of functions $f$ such that $f\circ \Phi$ has a
  $2\pi$-periodic flow near $\Phi^{-1}(c)$ (the primitive period may
  be any $2\pi/k$, where $k\in\NM^*$). This defines a sheaf over $B$.
  By Liouville-Arnold-Mineur, and since the fibre $\Phi^{-1}(c)$ is
  connected, the stalk over a regular value is isomorphic to $\ZM^2$.
  By Eliasson's normal form, this also holds near an elliptic critical
  value. Since no other type of critical point occur, our sheaf is
  just a flat bundle over $B$, and since $B$ is simply connected,
  there is no obstruction to the existence of a global section of the
  associated frame bundle, which is by definition a smooth map
  $g:B\fleche\RM^2$ which is a local diffeomorphism and such that
  $g\circ \Phi$ defines an effective Hamiltonian $\T^2$-action on $M$.
  
  Conversely, if $\Phi$ is of toric type, we know from
  proposition~\ref{prop:diffeo} that the fibres are connected and the
  image $B=\Phi(M)$ (and even $\rond{B}$) is connected and simply
  connected. Now proposition~\ref{prop:elliptic} tells us that no
  focus-focus singularities are present. Hence by
  proposition~\ref{prop:fibres} we have $B_r=\rond{B}$ and hence is
  connected and simply connected.
\end{demo}

\section{Moment polygons for semi-toric momentum maps}

In this section we come to our main point, defining moment polyhedrons
(here, polygons) for a particular class of almost-toric momentum maps,
roughly speaking those for which an $S^1$ action persists.

To be precise, what we shall call a polygon is a closed subset of
$\RM^2$ whose boundary is a continuous, piecewise linear curve with a
finite number of vertices in any compact. A convex polygon is
equivalently the convex hull of isolated points in $\RM²$. A polygon
is rational is the difference of the slopes of consecutive edges is
always rational.

We assume throughout that $\Phi$ is almost-toric (which, we recall,
requires $\Phi$ proper). By Liouville-Arnold-Mineur, The image of
$\Phi$ is naturally endowed with an integral affine structure with
boundary (which means that the boundary is a piecewise linear curve,
where linear means geodesic with respect to the affine structure): by
definition, affine charts are action variables, \emph{ie.} maps
$f:U\fleche \RM^2$, where $U$ is a small open subset of the image of
$\Phi$ and $f\circ \Phi$ generates a Hamiltonian $\T^2$-action (more
precisely, each component of $f\circ\Phi$ need have a $2\pi$-periodic
Hamiltonian flow.) This affine structure is integral because any two
such charts differ by the action of the integral affine group
$\textup{GA}(n,\ZM):=\RM^2\rtimes\textup{GL}(n,\ZM)$. Many more
details can be found in~\cite{duistermaat,zung-II,san-dijon} or
even~\cite{kontsevich-soibelman}.

Moreover by Eliasson's normal form at elliptic-elliptic singularities
the corners are convex and rational in any affine chart (we will show
below that the fibres are connected, hence by
proposition~\ref{prop:fibres} the boundary is exactly the set of
elliptic critical values).  This is just due to the fact that a germ
of convex sector near its summit is sent by a local diffeomorphism to
a germ of convex sector.  So in this sense the image of $\Phi$ is
always a kind of rational convex polygon (with focus-focus critical
values inside).  To have a \emph{true} polygon in $\RM^2$ we need to
find a natural projection of the universal cover of $B_r$ onto
$\RM^2$, respecting the affine structure.

\begin{defi}
  We say that an almost-toric $\Phi$ is \textbf{semi-toric} if there
  is a local diffeomorphism $f=(f^{(1)},f^{(2)})$ on the image of
  $\Phi$ such that $f^{(1)}\circ\Phi$ is a proper momentum map for an
  effective action of $S^1$.
\end{defi}
The terminology \emph{semi-toric} may be confusing, with the risk of
being mistaken for \emph{almost-toric}. A more precise phrase would be
``almost-toric with deficiency index one'' or ``almost-toric with
complexity one''~\cite{karshon-tolman}.  We shall keep semi-toric for
its shortness.

\begin{prop}
  \label{prop:diffeo2}
  In the definition above, $f$ is a diffeomorphism from the image of
  $\Phi$ into the image of $f\circ\Phi$.
\end{prop}
\begin{demo}
  The proof is the same as that of proposition~\ref{prop:diffeo},
  provided we show that the fibres of $f\circ\Phi$ are connected. But
  this is shown by theorem~\ref{theo:image} below.
\end{demo}

\begin{rema}
  The condition that the momentum map for the $S^1$ action is proper
  is very strong. In our situation this implies in many cases that $M$
  is compact (this is due to the presence of focus-focus singularities
  --- see corollaries~\ref{coro:number} and \ref{coro:compact}), and
  compact symplectic 4-manifolds with such an action are classified
  by~\cite{audin-topology} and~\cite{karshon-S1}. On the other hand
  many situations in classical mechanics when focus-focus
  singularities appear do have a global $S^1$ action but with
  non-proper momentum map~: a famous example is the spherical
  pendulum. Our results are still relevant to these cases when one can
  perform a preliminary reduction, or symplectic
  cutting~\cite{lerman-cut}, or more generally some integrable
  surgery~\cite{zung-II}, which isolates the interesting part of the
  manifold, making the induced $S^1$ momentum map proper.  Since we
  are more interested in the momentum map rather than in the
  symplectic manifold itself, this is a quite harmless operation.
\end{rema}

Assume now that $\Phi$ is semi-toric. We switch to the new momentum
map $f\circ\Phi$, which we call $\Phi$ again, and denote by $(J,H)$
its components. We define $J_{\textup{min}}$ (resp.
$J_{\textup{max}}$) to be the (possibly infinite) minimum (resp.
maximum) of $J$ on $M$.

\begin{theo}
  \label{theo:image}
  \begin{enumerate}
  \item The functions $H^+(x):=\max_{J^{-1}(x)}H$ and
    $H^-(x):=\min_{J^{-1}(x)}H$ are continuous;
  \item The image $B=\Phi(M)$ is the domain defined by
    \begin{equation}
      B=\{(x,y)\in\RM^2,\quad J_{\textup{min}}\leq x\leq J_{\textup{max}}
      \textrm{ and } H^-(x)\leq y\leq H^+(x)\}.
      \label{equ:image}
    \end{equation} 
    (Therefore $B$ is simply connected.)
  \item The fibres $\Phi^{-1}(c)$ are connected (and therefore the
    critical values of $\Phi$ are described as in
    proposition~\ref{prop:fibres}).
  \end{enumerate}
\end{theo}
\begin{demo}
  \emph{1}. By standard Morse theory, a discontinuity of $H^{+,-}$
  could only appear at a critical value of $J$. But since $B$ is
  closed this means that $\partial B$ would have a vertical segment in
  it. By hypothesis a segment of critical values can only correspond
  to a family of codimension 1 elliptic singularities of $\Phi$. Hence
  the preimage of this segment by $\Phi$ would be a locally maximal or
  minimal manifold for $J$, which is impossible (except for
  $x=J_{\textup{min}}$ or $x=J_{\textup{max}}$): by Morse-Bott theory
  (see \cite{atiyah-convex}) $J$ has a unique locally maximal (resp.
  minimal) manifold.

  \emph{2}.  Since $J$ is a proper momentum map for a Hamiltonian
  $S^1$ action on $M$, the fibre $J^{-1}(x)$ is compact and connected.
  Hence $H(J^{-1}(x))$ is compact and connected. Since by definition
  \[
  B= \bigsqcup_{x\in[J_{\textup{min}},J_{\textup{max}}]}\{x\}\times
  H(J^{-1}(x)),
  \]
  we have the description~\eqref{equ:image}. In particular $B$ is
  contractible to a line segment and hence is simply connected.

  Finally, to prove the connectedness statement \emph{3}, we still
  proceed similarly to \cite{atiyah-convex} (even if we are not in the
  toric case). For a regular value $x$ of $J$, $Z:=J^{-1}(x)$ is a
  smooth compact connected manifold. By the non-degeneracy hypothesis,
  $H_{\restr Z}$ is a Morse-Bott function with index $0$ or $2$. Hence
  the fibres of $H_{\restr Z}$ are connected. By continuity all fibres
  of $\Phi$ are connected.
\end{demo}
Notice that it is quite remarkable that semi-toric implies
connectedness of the fibres, as in the standard toric theorem where
\emph{both} Hamiltonians $H$ and $J$ needed to be periodic.

\begin{figure}[htbp]
  \begin{center}
    \input{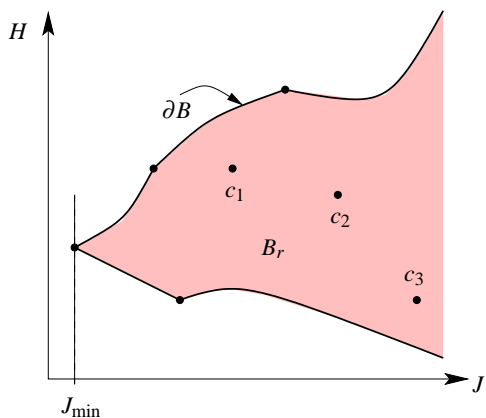}
    \caption{Image of $\Phi$}
    \label{fig:image}
  \end{center}
\end{figure}
\begin{coro}
  If $\Phi$ is semi-toric then $\Phi$ is of toric type if and only if
  it has no focus-focus singularity.
\end{coro}
\begin{demo}
  Combine the theorem with propositions~\ref{prop:fibres}
  and~\ref{prop:simply-connected}.
\end{demo}
\begin{rema}
  Some of the proofs above could be made even more natural by
  considering the symplectic reduction by $J$, and also the so-called
  symplectic cutting. If $x$ is a regular value of $J$, then we
  restrict $\Phi$ to the symplectic submanifold equal to
  $J^{-1}([x-\epsilon, x+\epsilon])$ with its boundary collapsed by
  the $S¹$ action. This manifold would by toric by
  proposition~\ref{prop:simply-connected}, and everything would follow
  from the standard toric theory. This would just require to state
  everything in the orbifold setting (which is probably not a big
  trouble in principle), since the action defined by $J$ is not
  necessarily free.
\end{rema}

\paragraph{Monodromy of focus-focus points. --- }
By the theorem, the local topological structure of the fibration by
$\Phi$ can be read off from the image $B$ (together with the
focus-focus critical values). It is true for the local symplectic
structure as well if one takes into account the integral affine
structure of $B$. As we said before, the affine structure on set $B_r$
of regular values of $\Phi$ comes from standard action variables; it
is extended on the boundary $\partial B$ using elliptic normal forms.
And its behaviour at focus-focus singularities is well understood.  In
particular one can compute the holonomy of this affine structure
around a focus-focus critical value $c_i$.  Recall from (for instance)
\cite{san-mono} that this holonomy (usually called the \emph{affine
  monodromy}) $\mu_A$ of $B_r$ is defined from a developing map as
follows. On the universal cover $\tilde{B_r}\flechedroite{\pi}B_r$ one
can define a global set of action variables (ie. a global affine
chart, or a developing map) $\tilde{f}:\tilde{B_r}\fleche\RM²$.  Let
$\gamma:[0,1]\fleche B_r$ be a loop starting at a point $c$ and
$\tilde{\gamma}$ a lift to $\tilde{B_r}$.  Then $\mu_A(\gamma,c)$ is
defined to be the element in $\textup{Aff}(2,\ZM)$ such that
\begin{equation}
  \tilde{f}(\tilde{\gamma}(0))=\mu_A(\gamma,c)\tilde{f}(\tilde{\gamma}(1))
  \label{equ:mono-affine}
\end{equation}
We know that the affine monodromy around a focus-focus critical value
$c_i$ has a unique line $\mathcal{L}$ of fixed points in $B$ (by line
we mean a geodesic of the affine structure). But given any
distribution of affine directions in $B$ we can associate a
1-dimensional vector space of locally Hamiltonian vector fields on
$M$: if $\beta$ is a closed 1-form on $B$ whose kernel gives the
affine directions, then we choose the symplectic dual of
$\Phi^*\beta$. In our case the smallest integral vector field
$\ham{1}$ corresponding to the direction of $\mathcal{L}$ is the
unique (up to sign) invariant Hamiltonian vector field generating an
$S^1$ action in a neighbourhood of the critical fibre.  More precisely,
suppose $c$ is a point close to $c_i$ and use such an
$\ham{1}=\ham{1}(c)$ to construct a basis $(\ham{1}(c),\ham{2}(c))$
corresponding to an integral affine basis $\mathcal{B}$ of $T_cB_r$;
next endow $\RM^2$ with the affine structure characterised by the
origin $f(c_1)$ and the basis $df(c).\mathcal{B}$. Then the affine
monodromy of an oriented loop $\gamma$ starting at $c$ and winding
once around $c_i$ has no translation component and its linear part is
equal to the matrix
\begin{equation}
  T^k:=\left(
    \begin{array}{cc}
      1 & 0\\ k & 1
    \end{array}
  \right)
  \label{equ:Tk}
\end{equation}
for some $k\in\NM$, which is the number of focus-focus critical points
in the critical fibre~\cite{matveev,zung-focus,cushman-san}.  The fact
that there is no translation component follows from the existence of a
symplectic potential in a neighbourhood of the critical fibre
$\Phi^{-1}(c_i)$ --- which in turn is due to the Lagrangian nature of
the critical fibre.

\begin{rema}
  While the holonomy of the affine structure thus determines the
  topology of the critical fibre, the semi-global symplectic
  classification of the fibration near $\Phi^{-1}(c_i)$ is a much
  harder issue. In this case it is non-trivial and given in the
  article~\cite{san-semi-global}.
\end{rema}

\paragraph{Structure of the image of $\Phi$. ---}
A developing map $\tilde{f}$ on $\tilde{B_r}$ can be uniquely extended
to the boundary $\pi^{-1}(\partial B)$.  Using the definition of the
affine structure on the boundary (given by the normal form of elliptic
singularities) we see that the boundary $\pi^{-1}(\partial B)$ is sent
by $\tilde{f}$ to a piecewise linear curve in $\RM^2$.  But the image
of $\tilde{f}$ is in general too non-injective to be of interest. In
our case, instead of going to the universal cover, it is easier to
make $B_r$ simply connected by suitable \emph{cuts}, and the image
obtained thereby becomes simple to interpret.

Let $\{c_i=(x_i,y_i),i=1,\dots,m_f\}\in\RM^2$ the set of focus-focus
critical values, ordered in such a way that $x_1\leq x_2\leq\cdots\leq
x_{m_f}$. For simplicity we have assumed that $m_f<\infty$, otherwise
just label $\{c_i\}$ with $i\in\ZM$, $\ZM^+$ or $\ZM^-$ and the rest
would go through. But we prove in corollary~\ref{coro:finite} below
that $m_f$ is actually always finite...

For each $i$ and for some $\epsilon\in\{-1,+1\}$ we define
$\mathcal{L}_i^{\epsilon}$ to be the vertical half line starting at
$c_i$ and going to $\epsilon\infty$:
$\mathcal{L}_i^{\epsilon}=\{(x_i,y), \epsilon y\geq\epsilon y_i\}$.
Given
$\vec\epsilon=(\epsilon_1,\dots,\epsilon_{m_f})\in\{-1,+1\}^{m_f}$, we
define the line segment $\ell_i:=B\cap \mathcal{L}_i^{\epsilon_i}$,
and
\begin{equation}
  \ell^{\vec\epsilon}=\bigcup_i\ell_i,
  \label{equ:cuts}
\end{equation}
where in addition we decorate each $\ell_i$ with the multiplicity
$\epsilon_i k_i$, where $k_i$ is the number of critical points in the
fibre $\Phi^{-1}(c_i)$.  More precisely, if several $c_i$'s have the
same $x_i$-coordinate, $\ell_i$ is the union of all corresponding
segments and we decide that each point $c$ in the
union~\eqref{equ:cuts} acquires the sum of the multiplicities
involved, which we denote by $k(c)$. A point with multiplicity zero is
omitted.

Let $\A^2_{\ZM}$ be the plane $\RM^2$ equipped with its standard
integral affine structure. The group of automorphisms of $\A^2_{\ZM}$
is the integral affine group
$\textup{Aff}(2,\ZM)=GL(2,\ZM)\ltimes\RM^2$. We denote by
$\mathcal{T}$ the subgroup of $\textup{Aff}(2,\ZM)$ which leaves a
vertical line (with orientation) invariant. In other words an element
of $\mathcal{T}$ is a composition of a vertical translation and an
element of $\{T^k,k\in\ZM\}\subset GL(2,\ZM)$.

\begin{theo}
  \label{theo:polygon}
  Given any $\vec\epsilon\in\{-1,+1\}^{m_f}$, there exists a
  homeomorphism $f$ from $B$ to $f(B)\in\RM^2$
  such that
  \begin{enumerate}
  \item $f_{\restr (B\setminus\ell^{\vec\epsilon})}$ is a
    diffeomorphism (into its image).
  \item $f_{\restr (B_r\setminus\ell^{\vec\epsilon})}$ is affine: it
    sends the integral affine structure of $B_r$ to the standard
    structure of $\A^2_{\ZM}$.
  \item $f$ preserves $J$: ie $f(x,y)=(x,f^{(2)}(x,y))$.
  \item $f_{\restr (B_r\setminus\ell^{\vec\epsilon})}$ extends to a
    smooth multi-valued map from $B_r$ to $\RM^2$ and for any
    $i=1,\dots,m_f$ and any $c\in\rond{\ell_i}$ then
    \begin{equation}
      \lim_{\substack{(x,y)\fleche c\\x<x_i}}df(x,y) =
      T^{k(c)}\lim_{\substack{(x,y)\fleche
          c\\x>x_i}}df(x,y),
      \label{equ:limit}
    \end{equation}
    where $k(c)$ is the multiplicity of $c$.
  \item The image of $f$ is a rational convex polygon.
  \end{enumerate}
  Such an $f$ is unique modulo a left composition by a transformation
  in $\mathcal{T}$.
\end{theo}
\begin{demo}
  We cannot show separately each point in the theorem. However we
  shall split the proof into several important steps.
  \paragraph{0.---}
  First of all, we use the description of the image of $\Phi$ given by
  theorem~\ref{theo:image} (and point \emph{4.} of
  proposition~\ref{prop:fibres}).  One can assume $m_f>0$. The case
  $m_f=0$ follows by proposition~\ref{prop:simply-connected} from the
  standard toric theory (and an argument like paragraph (2.---)
  below).
  
  \paragraph{1.---}
  For $i=0,\dots,m_f$, let $I_i$ be the open interval
  $(x_{i},x_{i+1})$ and (if $I_i\neq\emptyset$) $M_i=J^{-1}(I_i)$,
  where by convention $x_0=J_{\textup{min}}\in\{-\infty\}\cup\RM$ and
  $x_{m_f+1}=J_{\textup{max}}\in\RM\cup\{+\infty\}$.
  Each $M_i$ is an (open) symplectic manifold endowed with the
  momentum map $\Phi_{\restr M_i}$, and on which the set of regular values of
  $\Phi$ is connected and simply connected (by
  theorem~\ref{theo:image}). Moreover, the critical points of
  $\Phi_{\restr M_i}$ are non-degenerate and of elliptic type. Thus,
  as in proposition~\ref{prop:simply-connected}, we can define global
  action coordinates: there exists a smooth map $f_i:M_i\fleche\RM^2$
  which is a diffeomorphism into its image $B_i:=\Phi(M_i)$ and such
  that $f\circ\Phi$ is momentum map for a torus action on $M_i$.

  \paragraph{2.---}
  Actually, since $J_{\restr M_i}$ already defines an $S^1$ action,
  there exists an integer $p\neq 0$ such that $\ham{J}/p$ can be
  chosen to be the first element of an integral basis of the period
  lattice defining action variables. In other words, $f_i$ can be
  chosen of the form $f_i(x,y)=(x/p,f_i^{(2)}(x,y))$. But then one can
  see that the action of $J$ is effective if and only if $p=\pm 1$ (we
  leave this to the reader).  Therefore one can always chose
  $f_i(x,y)=(x,f_i^{(2)}(x,y))$.
  
  \paragraph{3.---}
  Assume first for simplicity that all $x_i$'s are different. Then
  $B\setminus\ell^{\vec\epsilon}$ is connected and simply connected.
  Since ${(f_0)}_{\restr B_0\cap B_r}$ is a section of the previously
  introduced sheaf of basic action variables on $B_r$, $f_0$ can be
  uniquely extended to a global section $f$ over
  $B\setminus\ell^{\vec\epsilon}$, and $J$ is always the first action
  variable.
  
  \paragraph{4.---}
  Remark that $B\setminus\ell^{\vec\epsilon}$ can be seen as a
  fundamental domain for the universal cover of $B_r$, and $f$ is a
  developing map for the affine structure. We look now at what happens
  at the gluing between $B_i$ and $B_{i+1}$ (fix $i=0$ for notational
  simplicity). Recall that in a neighbourhood of a focus-focus
  singularities there is a unique (up to a sign) Hamiltonian vector
  field $\ham{1}$ tangent to the fibres and whose flow is
  $2\pi$-periodic. And this vector field corresponds to a line through
  $c_1$ which is fixed by the affine monodromy (see paragraph ???). In
  our situation $\ham{1}$ must be $\pm\ham{J}$ and hence the fixed
  line $\mathcal{L}$ is the vertical line through $c_1$.  This implies
  that $f$ is continuous at $\ell_1\setminus\{c_1\}$, and the
  characterisation~\eqref{equ:limit} follows
  from~\eqref{equ:mono-affine}.
  \begin{figure}[htbp]
    \begin{center}
      \input{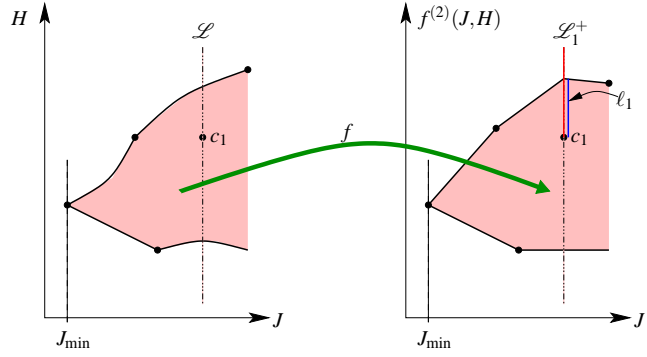}
      \caption{Definition of $f$ at $\ell_1$}
      \label{fig:cut}
    \end{center}
  \end{figure}
  
  We prove now that $f$ extends to a continuous map at $c_1$.  For
  this one can use the local normal form of~\cite{san-semi-global}.
  Since $f^{(2)}(J,H)$ is an action variable in $U\setminus\ell_1$,
  where $U$ is a neighbourhood of $c_1$, it follows from~\cite[remark
  3.2]{san-semi-global} that in coordinates $(\tilde{x},\tilde{y})$ of
  the form $\tilde{x}=x,\tilde{y}=\phy(x,y)$ for some function
  $\phy\in\Cinf(\RM^2,0)$,
  \[
  f^{(2)}(x,y)=\tilde{y}\ln\abs{\tilde{z}} - \tilde{x}\arg \tilde{z} +
  g(\tilde{x},\tilde{y}),
  \]
  where $(x,y)\in U\setminus\ell_1$,
  $\tilde{z}:=\tilde{x}+i\tilde{y}$, and $g$ is smooth at the origin.
  This shows that the function equal to $f^{(2)}$ in
  $U\setminus\ell_1$ and to $\tilde{y}\ln\abs{\tilde{z}} +
  g(0,\tilde{y})$ on $U\cap\ell_1$ is continuous in $U$.

  \paragraph{5.---}
  Notice that our construction of $f$ amounts to saying that $f_0$ on
  $B_0$ has been extended to $B_1$ by following paths in $B_r$ whose
  rule is to go only \emph{below} $c_1$ or \emph{above} $c_1$
  (depending on the sign of $\epsilon_1$). If several $x_i$'s are
  equal, one cannot necessarily find a path that goes only below some
  $c_i$ and above some others (in other words,
  $B\setminus\ell^{\vec\epsilon}$ is not necessarily connected). But
  we shall do the following: chose an arbitrary order $i_1,\dots,i_n$
  for the indices $i$ with the same value of $x_i$.  Then there is a
  unique (up to homotopy) path that connects $B_0$ and $B_1$ avoiding
  the $c_i$'s such that the whole picture is isotopic to a one where
  $x_{i_1}<x_{i_2}<\cdots<x_{i_n}$ and the path respects the above
  rule (see fig.~\ref{fig:plusieurs}).  Since the monodromy is
  Abelian~\cite{cushman-san}, the choice of the ordering does not
  affect the definition of $f$ in $B_1$, and the results follow as
  well.
  \begin{figure}[htbp]
    \begin{center}
      \input{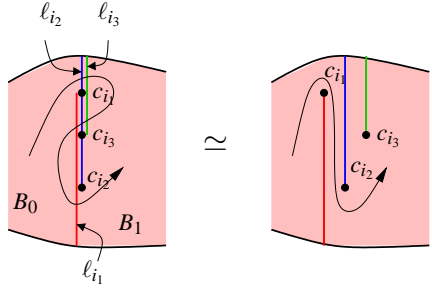}
      \caption{Extension of $f_0$ to $B_1$}
      \label{fig:plusieurs}
    \end{center}
  \end{figure}

  \paragraph{6.---}
  Since $f\circ\Phi$ is momentum map for a torus action on $M_i$, the
  boundary of $f\circ\Phi(M_i)$ corresponding to elliptic
  singularities is piecewise linear. Thus, the functions
  $[J_{\textup{min}},J_{\textup{max}}]\ni x\fleche
  f^{(2)}(x,H^\pm(x))$ (with the notations of
  theorem~\ref{theo:image}) are piecewise linear with rational slopes,
  and we have shown in paragraph (4.---) that they are continuous. It
  remains to show that the polygon $f\circ\Phi(M)$ is convex, which
  amounts to prove that $f^{(2)}(x,H^+(x))$ is convex and
  $f^{(2)}(x,H^-(x))$ is concave.  For this it suffices to look at the
  vertices. At elliptic-elliptic critical values, the result follows
  from the normal form. The other vertices that can appear are the
  points $(x_i,f^{(2)}(x_i,H^\pm(x_i)))$. Let us look for instance at
  the image of $v_1:=(x_1,H^+(x_1))$. Up to a change of sign for
  $f^{(2)}$, one can assume that $f(v_1)$ is still on the top boundary
  (which says that $f$ preserves the orientation). Let
  $\alpha=\lim_{\substack{x\fleche x_1\\x<x_1}}(H^+)'(x_1)$, ie.  the
  slope of the left-hand tangent to the boundary of $\Phi(M)$ at
  $v_1$, and $\beta$ the slope of the right-hand tangent. If $v_1$ is
  not the image of an elliptic-elliptic critical point then
  $\alpha=\beta$, otherwise $\beta<\alpha$ (the precise relation
  between $\alpha$ and $\beta$ is not needed here but will be given in
  section~\ref{sec:DH} below). Call $\alpha'$ and $\beta'$ the
  corresponding slopes for the new ``momentum map'' $f\circ\Phi$. (In
  other words they are the slopes of the edges of our moment polygon
  connecting at $f(v_1)$). Using~\eqref{equ:limit} we compute
  \[
  \beta'-\alpha' = \lim_{\substack{(x,y)\fleche
      v_1\\x>x_1}}\left(\deriv{f^{(2)}}{x}+\beta\deriv{f^{(2)}}{y}\right)
  - \lim_{\substack{(x,y)\fleche
      v_1\\x<x_1}}\left(\deriv{f^{(2)}}{x}+\alpha\deriv{f^{(2)}}{y}\right)
  \]
  \begin{equation}
    {}= -k(v_1) + (\beta-\alpha)\deriv{f^{(2)}}{y}(v_1).
    \label{equ:eq-slopes}
  \end{equation}
  Since $f$ is orientation preserving, one has
  $\partial{f^{(2)}}/\partial{y}>0$, hence
  \begin{equation}
    \beta'-\alpha'  \leq -k(v_1).
    \label{equ:ineq-slopes}
  \end{equation}
  Since $k(v_1)\geq 0$ (the cuts $\ell_i$ that can attain $v_1$ are
  only those that go up: for which $\epsilon_i=1$), the polygon is
  locally convex at the vertex $f(v_1)$ (or possibly flat if there is
  no cut and $\alpha=\beta$).
  
  Finally, if $v_1$ is an elliptic-elliptic vertex then
  $f\circ\Phi_{M_1}$ extends naturally to a smooth momentum map near
  $\Phi^{-1}(v_1)$ that gives local action coordinates. Hence we know
  as before that the slopes of the boundary of the local angular
  sector obtained by this momentum map are rational. This means that
  the last term in~\eqref{equ:eq-slopes} is rational, and
  $\beta'-\alpha'$ is thereby always rational.
  
  The other cases are handled in the same way, modulo only some sign
  changes.
\end{demo}

\begin{rema}
  As I learned afterwards, the use of such branch cuts was also crucial
  in Symington's work~\cite{symington-four}. They were called ``branch
  curves''; switching from upward to downward or vice-versa is a
  special case of her ``branch moves''.
\end{rema}
\section{The group of polygons}

Let $M$ be a symplectic 4-manifold equipped with a semi-toric momentum
map $\Phi=(J,H)$ with $m_f$ focus-focus critical fibres. For any
$\vec\epsilon\in\{-1,+1\}^{m_f}$ theorem~\ref{theo:polygon} gives an
equivalence class of rational convex polygons that we denote by
$\P_{\vec\epsilon}$, where the equivalence is given by the action of
transformations in $\mathcal{T}$. If one changes $\vec\epsilon$ the
class of $\P_{\vec\epsilon}$ modulo $\mathcal{T}$ might change. We
investigate here the relations between all these classes of polygons.

Given an affine vertical line $\mathcal{L}\subset\RM^2$ and an integer
$n\in\ZM$ we define a piecewise affine transformation
$t^n_{\mathcal{L}}$ of $\RM^2$ as follows: $\mathcal{L}$ splits
$\RM^2$ into two half-spaces. $t^n_{\mathcal{L}}$ acts as the identity
on the left half-space, and as the matrix $T^n$ (defined in
equation~\eqref{equ:Tk}) on the right one, for an origin of the affine
plane $\RM^2$ placed arbitrarily in $\mathcal{L}$ (recall that $T^n$
fixes $\mathcal{L}$).

We consider now the vertical lines $\mathcal{L}_i$ through the
focus-focus critical values $c_1,\dots,c_{m_f}$, and for any $\vec
n:=(n_1,\dots,n_{m_f})\in\ZM^{m_f}$ we construct the piecewise affine
transformation of $\RM^2$ $t_{\vec
  n}:=t^{n_1}_{\mathcal{L}_1}\circ\cdots\circ
t^{n_{m_f}}_{\mathcal{L}_{m_f}}$. This defines an Abelian action of
$\ZM^{m_f}$ on $\RM^2$.  Finally let $G=\{0,1\}^{m_f}$ (viewed as the
Abelian group $(\ZM/2\ZM)^{m_f}$) and let $\vec k=(k_1,\dots,k_{m_f})$
where $k_i\in\NM^*$ is the number of focus-focus critical points in
the fibre $\Phi^{-1}(c_i)$.

\begin{prop}
\label{prop:group}
  Let $G$ acts transitively on the set
  \[
  \P:=\{\P_{\vec\epsilon}, \vec\epsilon\in\{-1,+1\}^{m_f}\}
  \]
  by the formula
  \begin{equation}
    G\times\P\ni(\vec u,\P_{\vec\epsilon})\fleche
    \vec u\cdot\P_{\vec\epsilon} := \P_{(-2\vec
      u+1)\cdot\vec\epsilon}.
    \label{equ:action2}
  \end{equation}
  Then this action is given by the $t_{\vec n}$ transformations as
  follows:
  \begin{equation}
    \vec u\cdot\P_{\vec\epsilon} = t_{\vec
      u\cdot\vec\epsilon\cdot\vec
      k}\P_{\vec\epsilon},
    \label{equ:action1}
  \end{equation}
  The action is free if and only if the abscissae $x_i$'s of the
  focus-focus critical values are pairwise distinct.
\end{prop}
In the statement of the proposition, the dot $\cdot$ between
$\vec\epsilon$, $\vec u$ or $\vec k$ means pointwise multiplication in
$\ZM^{m_f}$, after the involved quantities $e_i\in\{-1,+1\}$ or
$u_i\in\{0,1\}$ are naturally injected in $\ZM$.  Notice that $\vec
u\fleche -2\vec u+1$ is just the standard group isomorphism between
$(\{0,1\}^{m_f},\oplus)$ and $(\{-1,+1\}^{m_f},\times)$, where
$\oplus$ is the addition modulo 2.

\begin{demo}
  The action $t_{\vec n}$ commutes with $\mathcal{T}$ and therefore
  induces an action on equivalences classes modulo $\mathcal{T}$.  It
  is then clear that both formulas \eqref{equ:action1} and
  \eqref{equ:action2} define actions of $G$ on some sets of
  equivalences classes modulo $\mathcal{T}$. So to prove that these
  actions coincide (and thereby that the result of \eqref{equ:action1}
  is indeed in $\P$) is suffices to look at generators of $G$.
  
  We use here the notations of the proof of
  theorem~\ref{theo:polygon}.  Selecting an element of the class of
  $\P_{\vec\epsilon} \mod \mathcal{T}$ amounts to fixing the starting
  local basis of action variables $f_0$ in $M_0$. Any other
  representative of that class can be obtained upon composing $f_0$ by
  a transformation in $\mathcal{T}$. So in what follows we fix $f_0$
  and by the notation $\P_{\vec\epsilon}$ we always mean the
  particular representative obtained from $f_0$ by the process of
  theorem~\ref{theo:polygon}.
  
  We assume here that the $x_i$'s are pairwise distinct. The general
  case follows, as before, by a splitting argument.
  
  Consider the action of $G$ given by equation~\eqref{equ:action2}: a
  ``1'' in the $ieth$ coefficient of $\vec u$ corresponds to a sign
  change in the $i$eth component of $\epsilon$, which flips the
  corresponding half line $\ell_i$ with respect to the point $c_i$.
  As a set of generators of $G$, we take the elements that have only
  one non-trivial coefficient.  Consider for instance the first one:
  $u=(1,0,\dots,0)$ and let it act on the polytope associated to the
  identity element $\vec 1$: $\vec u \cdot \P_{\vec 1} =
  \P_{\vec\epsilon}$, where $\epsilon=(-1,1,\dots,1)$.  Let $f_1$ and
  $\tilde{f}_1$ be the local action variables in $M_1$ obtained for
  $\P_{\vec 1}$ and $\P_{\vec\epsilon}$, respectively.  Let us fix for
  instance $y>y_1$ (recall that $c_1=(x_1,y_1)$ is a focus-focus
  critical value).  By~\eqref{equ:limit} one has for $\P_{\vec 1}$
  \[
  \lim_{\substack{(x,y)\fleche c\\x<x_i}}df_0(x,y) =
  T^{k_1}\lim_{\substack{(x,y)\fleche c\\x>x_i}}df_1(x,y),
  \]
  whereas for $\P_{\vec\epsilon}$ the formula reads
  \[
  \lim_{\substack{(x,y)\fleche c\\x<x_i}}df_0(x,y) =
  T^0\lim_{\substack{(x,y)\fleche c\\x>x_i}}d\tilde{f}_1(x,y),
  \]
  entailing
  \[
  \lim_{\substack{(x,y)\fleche c\\x>x_i}}d\tilde{f}_1(x,y) =
  T^{k_1}\lim_{\substack{(x,y)\fleche c\\x>x_i}}df_1(x,y),
  \]
  and therefore, since in $M_1$ $f_1$ and $\tilde{f}_1$ must differ
  only by an element of $\mathcal{T}$,
  \begin{equation}
    \label{equ:swap}
    \tilde{f}_1 = T^{k_1} \circ f_1.
  \end{equation}
  Now for $i>1$ the half lines $\ell_i$ are identical for
  $\P_{\vec\epsilon}$ and $\P_{\vec 1}$; this means that both $f_1$
  and $\tilde{f}_1$ are extended further in the same way, ensuring
  that for all $i>1$, $\tilde{f}_i = T^{k_1} \circ f_i$. This in turn
  says that the polytopes are precisely related by the formula
  $\P_{\vec\epsilon} = t_{\vec u \cdot \vec k} \P_{\vec 1}$. Doing
  this for all generators $\vec u$ we have proved that for all $\vec
  u\in G$,
  \[
  \P_{(-2\vec u+1)\cdot\vec 1} = t_{\vec u \cdot \vec k} \P_{\vec 1}.
  \]
  
  We conclude for a general $\vec \epsilon$ by the following
  elementary chasing around: let
  $\phy:(\{-1,+1\}^{m_f},\times)\fleche(\{0,1\}^{m_f},\oplus)$ be the
  isomorphism used in the statement of the proposition:
  $\phy^{-1}(\vec u)=-2\vec u+ 1$.  Thus one can write
  $\P_{\vec\epsilon}=\P_{\vec\epsilon\cdot \vec
    1}=t_{\phy(\vec\epsilon)\cdot\vec k}\P_{\vec 1}$.  Therefore
  $\P_{\phy^{-1}(\vec u)\cdot\vec\epsilon}=\P_{\phy^{-1}(\vec
    u)\cdot\vec\epsilon\cdot\vec 1} = \P_{\phy^{-1}(\vec
    u\oplus\phy(\vec\epsilon))\cdot\vec 1}=t_{\vec
    u\oplus\phy(\vec\epsilon)\cdot\vec k - \phy(\vec\epsilon)\cdot\vec
    k}\P_{\vec\epsilon}$. Now it is straightforward to check that
  $(\vec u\oplus\phy(\vec\epsilon)-\phy(\vec\epsilon))=\vec
  u\cdot\vec\epsilon$.  This shows that the right hand sides of
  \eqref{equ:action1} and \eqref{equ:action2} are indeed equal.
  
  The transitivity of the action is ensured by
  $\P_{\vec\epsilon}=\phy(\vec\epsilon)\cdot \P_{\vec 1}$. Finally,
  the subgroup of affine transformations generated by
  $t_{\mathcal{L}}$ acts freely on the set of all non-vertical
  segments starting on the right of $\mathcal{L}$. Applying this fact
  to the edges of the polygons $\P_{\vec\epsilon}$ one sees that the
  action of $G$ on $\P$ is free provided the $x_i$'s are distinct. Now
  suppose $x_{i}=x_{i+1}=\cdots=x_{i+j}$. Then the order in which we
  consider $c_i,\dots,c_{i+j}$ is irrelevant, and the corresponding
  permutation group in $j+1$ elements acts trivially on $\P$. In
  particular the action of $G$ is not free.
\end{demo}
\begin{rema}
  Let $\tilde{\P}$ be the set of all possibles polygons obtained for a
  given semi-toric momentum map $(J,H)$. As remarked in the proof,
  fixing a starting set of action variables $f_0$ gives a way of
  selecting a representative in each class $\P_{\vec\epsilon}$. This
  says that $\tilde{\P}$ is in bijection with $\P\times\mathcal{T}$,
  acquiring thereby a natural group structure, where the identity
  element if the representative of the class $\P_{\vec 1}$. In other
  words one has a short exact sequence
  \[
  0\fleche\mathcal{T}\fleche\tilde{\P}\fleche\P\fleche 0,
  \]
  which has a cross section given by the choice of $f_0$. If all the
  $x_i$'s are distinct then $\tilde{\P}$ is isomorphic to
  $G\times\mathcal{T}$.
\end{rema}

\section{Duistermaat-Heckman measures}%

\subsection{The $S¹$ action}
\label{sec:DH}
The polygons introduced in theorem~\ref{theo:polygon} are a very
efficient tool for recovering the various invariants attached to the
momentum map $\Phi$, and in particular to the effective $S¹$ action
defined by $J$.

We consider here the standard Duistermaat-Heckman measure $\mu_J$ for
the Hamiltonian $J$. Recall that by definition
$\mu_J([a,b])=\textrm{vol}(J^{-1}([a,b]))$, where $\textrm{vol}$ means
the symplectic (or Liouville) volume in $M$. It is known (see
\cite{duist-heckman}) that
\[
\mu_J:=\rho_J(x)\frac{\abs{dx}}{2\pi},
\]
where the density $\rho_J(x)$ (sometimes called the
Duistermaat-Heckman function) is a continuous function, equal to the
symplectic volume of the reduced orbifold $J^{-1}(x)/S¹$.
\begin{prop}
  \label{prop:DH}
  Given any $\vec\epsilon\in\{-1,1\}^{m_f}$ and any polygon $P$ in
  $\P_{\vec\epsilon}$, $\rho_J(x)$ is equal to the length of the
  vertical segment, intersection of the vertical line through $x$ and
  the (filled) polygon $P$. Hence $\rho_J(x)$ is piecewise linear.
\end{prop}
\begin{demo}
  Of course the fact that $\rho_J(x)$ is piecewise linear also follows
  from the theorem of Duistermaat and Heckman. It comes very easily
  here because we are in an integrable situation. Namely let $f$ be
  the homeomorphism given by theorem~\ref{theo:polygon}. Then in each
  ``cell'' $M_i$, $\tilde{\Phi}:=f(J,H)=(J,f^{(2)}(J,H))$ is a set of
  smooth action variables. If follows from Liouville-Arnold-Mineur
  theorem that the Duistermaat-Heckman measure $\mu_{\tilde{\Phi}}$ on
  $\RM²$ associated to $\tilde{\Phi}$ has density 1 over
  $\abs{dx\wedge dy}/(2\pi)²$. Integrating in the vertical direction
  one finds the result.
\end{demo}
\begin{rema}
  This shows that the lengths of the vertical segments of the polygons
  in $\P$ don't depend either on $\vec\epsilon$ or on the particular
  representative. This, of course, can also be checked directly from
  the definition of these polygon (the action of $\mathcal{T}$ does
  not change vertical lengths).
\end{rema}

We calculate now $\rho_J(x)$ in terms of the generalised moment
polygons of theorem~\ref{theo:polygon}. Let
$\vec\epsilon\in\{-1,+1\}^{m_f}$ and let $f=f_{\vec\epsilon}$ be the
homeomorphism given by the theorem. As before, $c_j$'s are the
focus-focus critical values and $k_j$ is the number of focus-focus
point in the fibre above $c_j$.

If $c$ is a critical value of maximal corank of $\Phi$, then
$\Phi^{-1}(c)$ is either of focus-focus point or an elliptic-elliptic
point. In the latter case we call $c$ a ``top vertex'' if it lies in
the graph of $H^+$ and a ``bottom vertex'' if it lies in the graph of
$H^-$ (in the terminology of theorem~\ref{theo:image}). At such a
critical point $J$ can be written in suitable symplectic coordinates
under the form $J=a(x²+\xi²)/2+b(y²+\eta²)/2$ for integer $a$, $b$
which are called \emph{isotropy weights} of the $S¹$ action defined by
$J$ \cite{duist-heckman,karshon-S1}.

\begin{theo}
  \label{theo:DH}
  If $\alpha^+(x)$ (resp. $\alpha^-(x)$) denotes the slope of the top
  (resp. bottom) boundary of the polygon $f\circ\Phi(M)$, then the
  derivative of the Duistermaat-Heckman function is
  \begin{equation}
    \rho'_J(x)=\alpha^+(x)-\alpha^-(x)
    \label{equ:derivative}
  \end{equation}
  and is locally constant on
  $J(M)\setminus\{\pi_x(f(\Sigma_0(\Phi)))\}\in\RM$, where
  $\Sigma_0(\Phi)$ is the set of critical values of $\Phi$ of maximal
  corank and $\pi_x$ is the projection $(x,y)\fleche x$. If
  $(x,y)\in\Sigma_0(\Phi)$ then
  \begin{equation}
    \rho'_J(x+0)-\rho'_J(x-0)= -\sum_j k_j - e^+ - e^-,
    \label{equ:DH}
  \end{equation}
  where the sum runs over the set of all indices $j$ such that
  $\pi_x(c_j)=x$, and $e^+$ (respectively $e^-$) is non-zero if and
  only if an elliptic top vertex (resp. a bottom vertex) projects down
  onto $x$. If this occurs then
  \[
  e^\pm = - \frac{1}{a^\pm b^\pm} \geq 0,
  \]
  where $a^\pm$, $b^\pm$ are the isotropy weights for the $S¹$ action
  at the corresponding vertices.
\end{theo}
\begin{demo}
  The first point is obvious in view of proposition~\ref{prop:DH}.
  Notice that in general the discontinuities of $\rho'_J$ occur at the
  singularities of $J$. Here these singularities (except possibly for
  the maxima and minima of $J$) are exactly critical values of maximal
  corank of $\Phi$.
  
  The second point is just a small refinement of
  formula~\eqref{equ:eq-slopes}. This formula says that
  \[
  \rho'_J(x+0)-\rho'_J(x-0)= -\sum_j k_j + (r^+(x)-r^-(x))
  \]
  where, as explained at the end of the proof of
  theorem~\ref{theo:polygon}, $r^\pm(x)$ is computed as follows. The
  item \emph{4.} of theorem~\ref{theo:polygon} says that in a small
  neighbourhood of the point in the boundary $\nu^\pm:=(x,H^\pm(x))$,
  $f$ can be smoothly extended (either from the region $\leq x$ or
  $\geq x$) to a smooth map $\tilde{f}^\pm$ such that
  $\tilde{f}^\pm\circ\Phi$ is a toric momentum map near
  $\Phi^{-1}(\nu^\pm)$. Then the local image of
  $\tilde{f}^\pm\circ\Phi$ is a angular sector and $r^\pm(x)$ is the
  difference between the slopes of the right-hand and left-hand edges
  at the vertex $f(\nu^\pm)$ of this sector. It does not depend on the
  way $f$ was extended since it is invariant by a transformation in
  $\mathcal{T}$. Precisely, there is a matrix
  $A^\pm=\left(\begin{array}{cc} a^\pm & b^\pm\\c^\pm & d^\pm
    \end{array}\right)\in SL(2,\ZM)$ and canonical coordinates
  $(x,y,\xi,\eta)$ near the elliptic-elliptic point
  $\Phi^{-1}(\nu^\pm)$ such that
  $\tilde{f}^\pm\circ\Phi=A^\pm\circ(\frac{x²+\xi²}{2},\frac{y²+\eta²}{2})$.
  In particular
  $J=a^\pm(\frac{x²+\xi²}{2})+b^\pm(\frac{y²+\eta²}{2})$.  If $x$ is
  not an extremal value for $J$, $a^\pm$ and $b^\pm$ do not vanish and
  have different signs; for the top vertex $\nu^+(x)$ one must have
  $a^+<0$.  Then
  \[
  r^+(x) = \frac{d^+}{b^+}-\frac{c^+}{a^+}=\frac{1}{a^+b^+}.
  \]
  At the bottom vertex $\nu^-(x)$ the coefficient $a^-$ is positive,
  and $r^-(x)=-\frac{1}{a^-b^-}$.
\end{demo}
\begin{rema}
  Nothing in this theorem is essentially new, apart from the proof
  (and maybe also the fact that $M$ is not necessarily compact).
  Compared to the usual theory, our proof follows very easily and
  elementarily from our moment polygons. For general Hamiltonian torus
  actions on compact symplectic manifolds, a formula analogous to
  \eqref{equ:DH} follows from the Duistermaat-Heckman formula (or the
  localisation formula of Atiyah-Bott-Berline-Vergne) for the Fourier
  transform of $\mu_J$, and a Fourier inversion argument as in
  \cite{guillemin-lerman-sternberg1} (see also \cite{karshon-S1}). The
  main difference with our formula is that we separate the
  contribution of focus-focus points from elliptic-elliptic points,
  which of course is not possible in the context of a general $S¹$
  action. This again is not really new since the link between the
  monodromy and Duistermaat-Heckman's theory was recently pointed out
  by Nguyên Tiên Zung in \cite{zung-another}.  However Zung's
  construction was a local one using integrable surgery, whereas we
  express it in a global situation.
\end{rema}

This theorem (together with theorem~\ref{theo:polygon}) has some easy
corollaries of topological nature.
\begin{coro}
  \label{coro:bounded}
  If a symplectic manifold $M$ admits a semi-toric momentum map
  $(J,H)$ with at least one critical value of maximal corank
  ($dJ(m)=dH(m)=0$) then $J$ is bounded from below or from above.
\end{coro}
\begin{demo}
  By the theorem \ref{theo:DH}, the strict inequality
  $\rho'_J(x+0)-\rho'_J(x-0)<0$ holds at at least one point.  Hence
  there is a point $x_0$ for which $\rho'_J(x_0)\neq 0$.  Suppose for
  instance $\rho'_J(x_0)<0$. Then the length of the interval
  $f^{(2)}(x,J^{-1}(x))$ (or the Duistermaat-Heckman measure at $x$)
  is bounded from above by $\textup{const}+x\rho'_J(x)$ and hence by
  convexity of the polygon must vanish for a finite value of $x>x_0$.
  The point for which it vanishes has to be the maximal value of $J$.
\end{demo}
\begin{coro}
  \label{coro:number}
  Let a symplectic manifold $M$ admit a semi-toric momentum map
  $\Phi=(J,H)$ such that $J$ is bounded from below with minimal value
  $J_{\textup{min}}$. If $\Phi$ has more than
  $\rho'_J(J_{\textup{min}}+0)$ focus-focus points (counted with
  multiplicity) then $M$ is compact.
\end{coro}
\begin{demo}
  By the theorem \ref{theo:DH} if $x$ is greater than the maximum of
  the abscissae of the focus-focus critical values then
  $\rho'_J(x)<0$, and we conclude as above that $J$ has a finite
  maximal value. Hence $M$ is compact by properness of $J$.
\end{demo}
\begin{rema}
  In case of a compact $M$, one can write an explicit upper bound for
  the symplectic volume of $M$ (the area of the polygon: see next
  section) in terms of $\rho'_J(J_{\textup{min}}+0)$, the symplectic
  volume of $J^{-1}(J_{\textup{min}})$ (which may be zero), and the
  abscissae and multiplicities of all focus-focus critical values. We
  leave this to the reader.
\end{rema}
\begin{coro}
  \label{coro:compact}
  If $M$ admits a semi-toric momentum map $\Phi=(J,H)$ with $m_f\geq
  2$ focus-focus critical fibres and such that $J$ has a unique
  minimum (or maximum) then $M$ is compact.
\end{coro}
\begin{demo}
  If $J$ has a unique minimum its image under $\Phi$ is an
  elliptic-elliptic corner of any associated moment polygon, open in
  the direction $y\geq 0$. But the edges of an elliptic-elliptic
  corner are directed along integral vectors $(a,c)$ and $(b,d)$ such
  that $\left(
    \begin{array}{cc}
      a & b\\c & d
    \end{array}
  \right)\in SL(2,\ZM)$. Hence
  $\rho'_J(J_{\textup{min}}+0)=d/b-c/a=1/ab\leq 1$ and the result
  follows from the corollary~\ref{coro:number} above.
\end{demo}
\begin{rema}
  In contrast with the hypothesis of this corollary,
  $\rho'_J(J_{\textup{min}}+0)$ can take any integral value if $J$ has
  a non-trivial submanifold of minima (which means that the moment
  polygons have a vertical edge at $J_{\textup{min}}$).
\end{rema}
\begin{coro}
  \label{coro:finite}
  If $M$ admits a semi-toric momentum map then the number $m_f$ of
  focus-focus critical fibres is finite.
\end{coro}
\begin{demo}
  If $m_f>0$ then by corollary~\ref{coro:bounded} $J$ is semi-bounded
  (say for instance from below). Then by corollary~\ref{coro:number}
  if $m_f$ is very large $M$ must be compact. Since the $c_i$'s are
  isolated $m_f$ has to be finite.
\end{demo}

\begin{rema}
  If one knows the value of $\rho_J$ at some point, then
  theorem~\ref{theo:DH} gives all one needs to reconstruct $\rho_J$ by
  integration. Such a formula will be given below for the
  ``generalised'' $S¹$ actions.
\end{rema}

\subsection{The generalised $S¹$ actions}

The construction of the polygons in theorem~\ref{theo:polygon} leads
naturally to considering another type of Duistermaat-Heckman measure,
namely the one associated with horizontal slices of the polygons. In
other words, given an $\vec\epsilon\in\{-1,1\}^{m_f}$ and a map $f$ as
provided by the mentioned theorem, we consider the push-forward of the
Liouville measure by the second ``action'' variable $K:=f^{(2)}(J,H)$.
$K$ is continuous, but it is not smooth along the vertical lines
through focus-focus points.  Where it is smooth, $K$ does define a
Hamiltonian $S¹$-action. The only problem for defining the
Duistermaat-Heckman measure $\mu_K$ is that $K$ is not assumed to be
proper. Therefore in what follows we either assume $M$ to be compact
or we restrict $M$ to the compact symplectic manifold with boundary
$J^{-1}([a,b])$, for some bounded interval $[a,b]$.

Then the Duistermaat-Heckman function $\rho_K$ such that
$\mu_K=\rho_K(y)\abs{dy}$ can be described exactly as we did for
$\rho_J$.  In particular $\rho_K(y)$ is the length of the horizontal
slice of the moment polygon with ordinate $y$. However this definition
is not so easy to use here since the order of the vertical projections
of the polygon vertices--- and hence $\rho_K(y)$ --- strongly depends
on the $\vec\epsilon$ chosen to construct it. It is more adequate to
express $\rho_K$ as much as possible in terms of the $J$-data.

For this purpose we slightly change the notation by letting
$x_0<x_1<\cdots<x_N$ be the absciss\ae{} of all critical values of
rank zero of $\Phi$ (including focus-focus and elliptic-elliptic
points).  Let $y_i^+$ (resp. $y_i^-$) be the ordinate of the
intersection of the vertical line at $x_i$ with the top (resp. bottom)
boundary of the polygon.  Finally for $i\in[0..N-1]$ let
$\alpha_i^\pm=\frac{y_{i+1}^\pm-y_{i}^\pm}{x_{i+1}-x_{i}}$ be the
slope of the corresponding (top or bottom) edge of the polygon.
Contrary to $x_i$, $\alpha^\pm_i$ and $y^\pm_i$ depend on
$\vec\epsilon$. One has~:
\[
y^\pm_{i}=y^\pm_0+\sum_{j=0}^{i-1} h_i\alpha^\pm_i, \textrm{ where }
h_i=(x_{i+1}-x_i),
\]
and $\alpha_{i+1}^\pm-\alpha_i^\pm$ is given by theorem~\ref{theo:DH}
in terms of $\vec\epsilon$, the monodromy indices, and fixed point
data of $J$.
\begin{figure}[htbp]
  \begin{center}
    \input{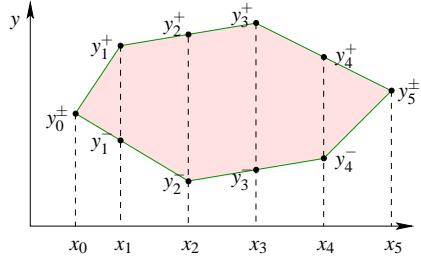}
    \caption{Notations for cutting the polygon}
    \label{fig:cutting}
  \end{center}
\end{figure}

We need some non-standard conventions in order to state the following
theorem. If $a,b$ are real numbers, we denote by $\intervalle{a,b}$
the interval $[\min(a,b),\max(a,b)]$. If $I$ is an interval, $\chi_I$
designates the characteristic function of $I$. If $I$ is a point, by
convention $\chi_I=0$, and for any number $\beta\in\RM\cup\{\infty\}$,
$\beta\chi_I=0$.

\begin{theo}
  With the notation defined above, the Duistermaat-Heckman function
  $\rho_K$ is the continuous, piecewise linear function given by the
  following formula~:
  \begin{equation}
    \label{equ:DH-function}
    \rho_K(y) = \sum_{i=0}^{N-1}
    \left(\frac{1}{\abs{\alpha^-_i}}(y-y_i^-)\chi_{\intervalle{y_i^-,y_{i+1}^-}}
      + h_i \chi_{\intervalle{y_{i+1}^-,y_i^+}} +
      \frac{1}{\abs{\alpha^+_i}}(y_{i+1}^+-y)\chi_{\intervalle{y_i^+,y_{i+1}^+}}\right).
  \end{equation}
  In particular the derivative of $\rho_K$ is the piecewise constant
  function given by
  \[
  \rho'_K(y) = \sum_{i=0}^{N-1}
  \left(\frac{1}{\abs{\alpha^-_i}}\chi_{\intervalle{y_i^-,y_{i+1}^-}}
    -
    \frac{1}{\abs{\alpha^+_i}}\chi_{\intervalle{y_i^+,y_{i+1}^+}}\right).
  \]
\end{theo}
\begin{demo}
  The term in the sum for a fixed $i$ corresponds to the calculation
  of $\rho_K$ restricted to the elementary cell
  $J^{-1}([x_i,x_{i+1}])$, which is a simple exercise.
\end{demo}

\section{Examples}%

\subsection{Coupled angular momenta on $S²\times S²$}
\label{sec:s2s2}
The first example that motivated this paper (with some others to
come), and which I still think is of primary interest, has been
described first by \sadov{} and \zhilin{} in
\cite{sadovski-zhilinski}. It is the problem of two coupled angular
momenta, describing for instance a so-called ``spin-orbit coupling''.
The momentum map on $S²\times S²$ depends on an additional parameter
$t$ as follows~: $\Phi_t=(J,H_t)$, where $J=N_z+S_z$ and
\[
H_t = \frac{1-t}{|\mathbf{S}|}S_z +
\frac{t}{|\mathbf{N}||\mathbf{S}|}\pscal{\mathbf{N}}{\mathbf{S}},
\quad 0\leq t\leq 1.
\]
We have denoted by $\mathbf{S}=(S_x,S_y,S_z)$ et
$\mathbf{N}=(N_x,N_y,N_z)$ the angular momentum variables on each $S²$
factor. In other words these spheres are standard symplectic spheres
but with radius $|\mathbf{S}|$ for the first one and $|\mathbf{N}|$
for the second one.

\begin{figure}[htbp]
  \centering
  \includegraphics{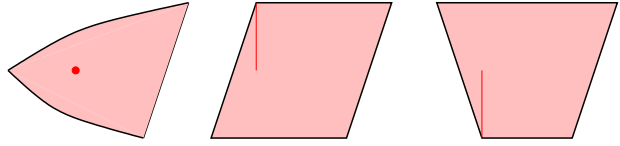}
  \caption{Image of the momentum map (left) and generalised polygons
    for the coupled angular momenta.}
  \label{fig:S2-S2}
\end{figure}

Then one can show that $\Phi_t$ is semi-toric except for two values of
$t$, and \emph{not} of toric type for $t$ in a bounded open interval
containing $1/2$, where $\Phi_t$ has a focus-focus critical point. For
$t$ around $1/2$ the image of the momentum map and the two generalised
polygons are depicted in the figure~\ref{fig:S2-S2}.  We don't show
the details here because they are partly computed
in~\cite{sadovski-zhilinski} and go along the same lines as the next
example.

\subsection{Coupled spin and oscillator on $S²\times\RM²$}
\label{sec:spin}
Using the previous example by \sadov{} and \zhilin, one can
construct an example on $S²\times\RM^2$ with one focus-focus
singularity, just by linearising one of the spheres at a pole. In
addition to being interesting by its computational simplicity, it
provides an example of a non-compact manifold that shows that
corollary~\ref{coro:compact} is optimal.

On $S²$ one has a natural Hamiltonian $S¹$ action whose Hamiltonian is
the ``vertical coordinate'' $z$, where we embed $S²$ in $\RM³$ as
$\{x²+y²+z²=1\}$. The sign of the ``standard'' symplectic form on $S²$
is chosen such that the flow turns around the vertical axis in the
direct sense (counterclockwise). The total symplectic volume is chosen
such that the flow of $z$ is $2\pi$-periodic.

On $\RM^2=\{(u,v)\}$ with canonical symplectic form our standard $S¹$
action is the harmonic oscillator $N:=(u²+v²)/2$ with $2\pi$-periodic
flow.

On $M=S²\times\RM²$ we define an $S^1$ action by the Hamiltonian
\[
J:=N+z.
\]

Using the embedding of $S²$ in $\RM³$, define the orthogonal projector
$\pi_z$ from $S²$ onto $\RM²$ viewed as the $z=0$ hyperplane. Let
$(m,p)\in S²\times\RM²$. Then under the flow of $J$ the points $m$ and
$p$ are moving along the flows of $z$ and $N$, respectively, with the
same angular velocity. Therefore the scalar product
$\pscal{\pi_z(m)}{p}$ is constant. That is,
\[
K:=(m,p)\fleche\pscal{\pi_z(m)}{p} = ux+vy
\]
commutes with $J$: $\{K,J\}=0$. Now we define
\[
H_t:=(1-2t)(N-z) + tK, \quad \textrm{ and } \quad \Phi_t:=(J,H_t).
\]

\begin{figure}[htbp]
  \begin{center}
    \includegraphics[width=7cm]{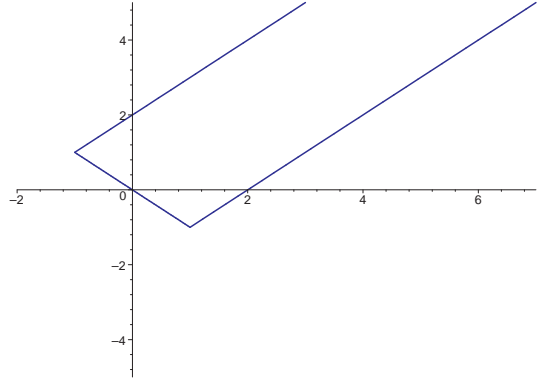}
    \caption{The standard moment polytope at $t=0$ for example
      \ref{sec:spin}.}
    \label{fig:moment_0}
  \end{center}
\end{figure}

When $t=0$, $\Phi_0=(N+z,N-z)$ defines an effective $\T²$ action and
hence is \emph{toric}. The moment polygon is depicted in
Fig.~\ref{fig:moment_0}. Notice that $\Phi_0$ is affinely equivalent
to the momentum map $(z,N)$ in which the variables are ``separated'',
or ``uncoupled''. Physically it describes a classical spin and a
harmonic oscillator. Hence the name we gave to this example (but it
probably deserves a better one). The particular linear scaling ($1-2t$
and $t$) is not important; it is just chosen in such a way that the
spectrum of the linearised Hamiltonian at the focus-focus point is
very simple.
\begin{prop}
  \begin{enumerate}
  \item For $t\in\RM\setminus\{1/3,1\}$ the momentum map $\Phi_t$ is
    semi-toric;
  \item for $t<1/3$ and $t>1$ the momentum map $\Phi_t$ is actually of
    toric type (in the sense of definition~\ref{defi:toric});
  \item for $t\in(1/3,1)$ the momentum map $\Phi_t$ is semi-toric with
    one simple focus-focus point;
  \item for $t\in\{1/3,1\}$ the momentum map $\Phi_t$ has a degenerate
    singularity (and hence is not almost-toric).
  \end{enumerate}
\end{prop}
If one needs only an example with one focus-focus point, the simplest
of course is to take $t=1/2$ or $\Phi=(J,K)$.
\begin{figure}[htbp]
  \begin{center}
    \emph{Image of the momentum map:}
    
    \includegraphics[width=2.2cm]{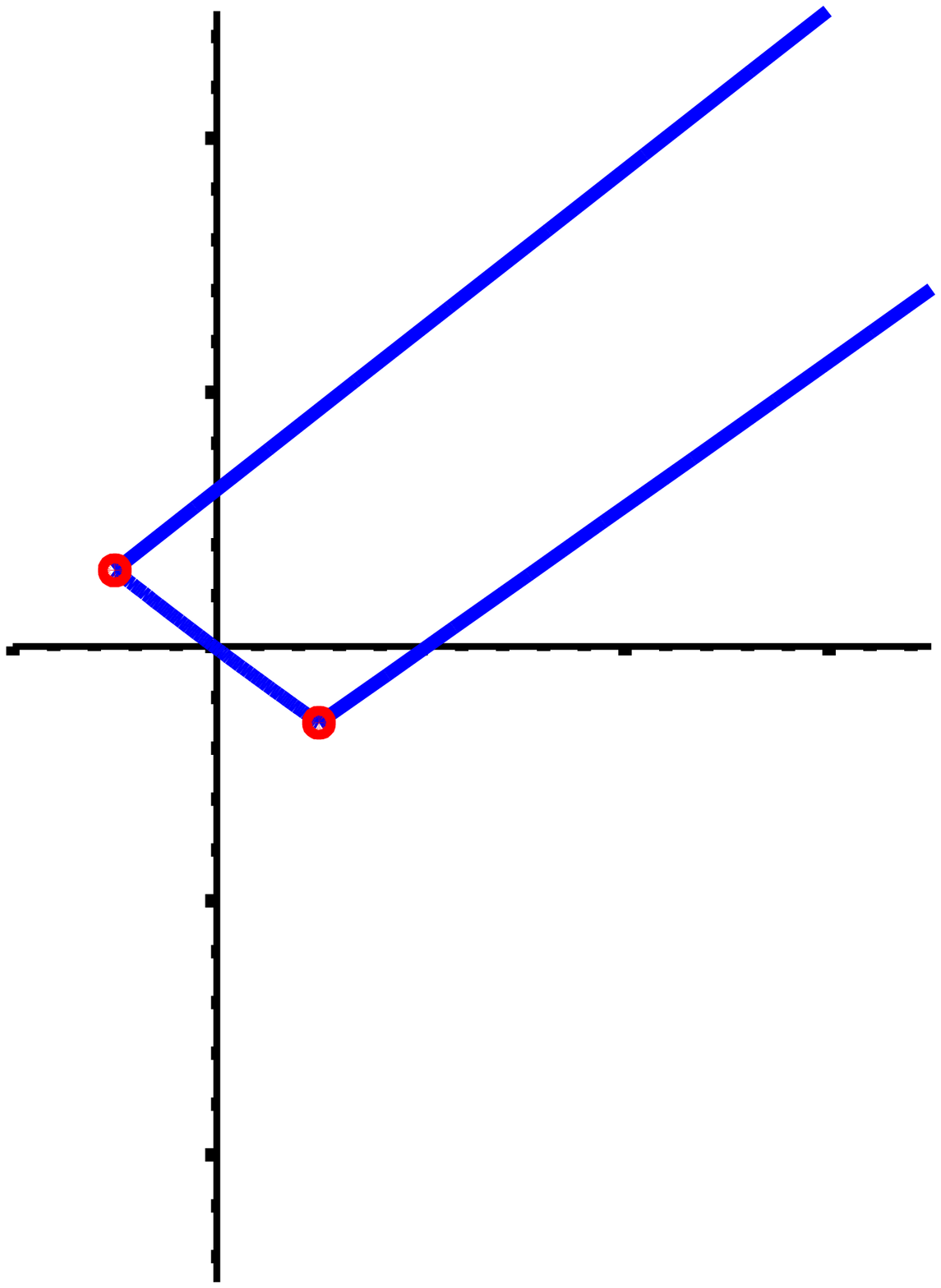}
    \includegraphics[width=2.2cm]{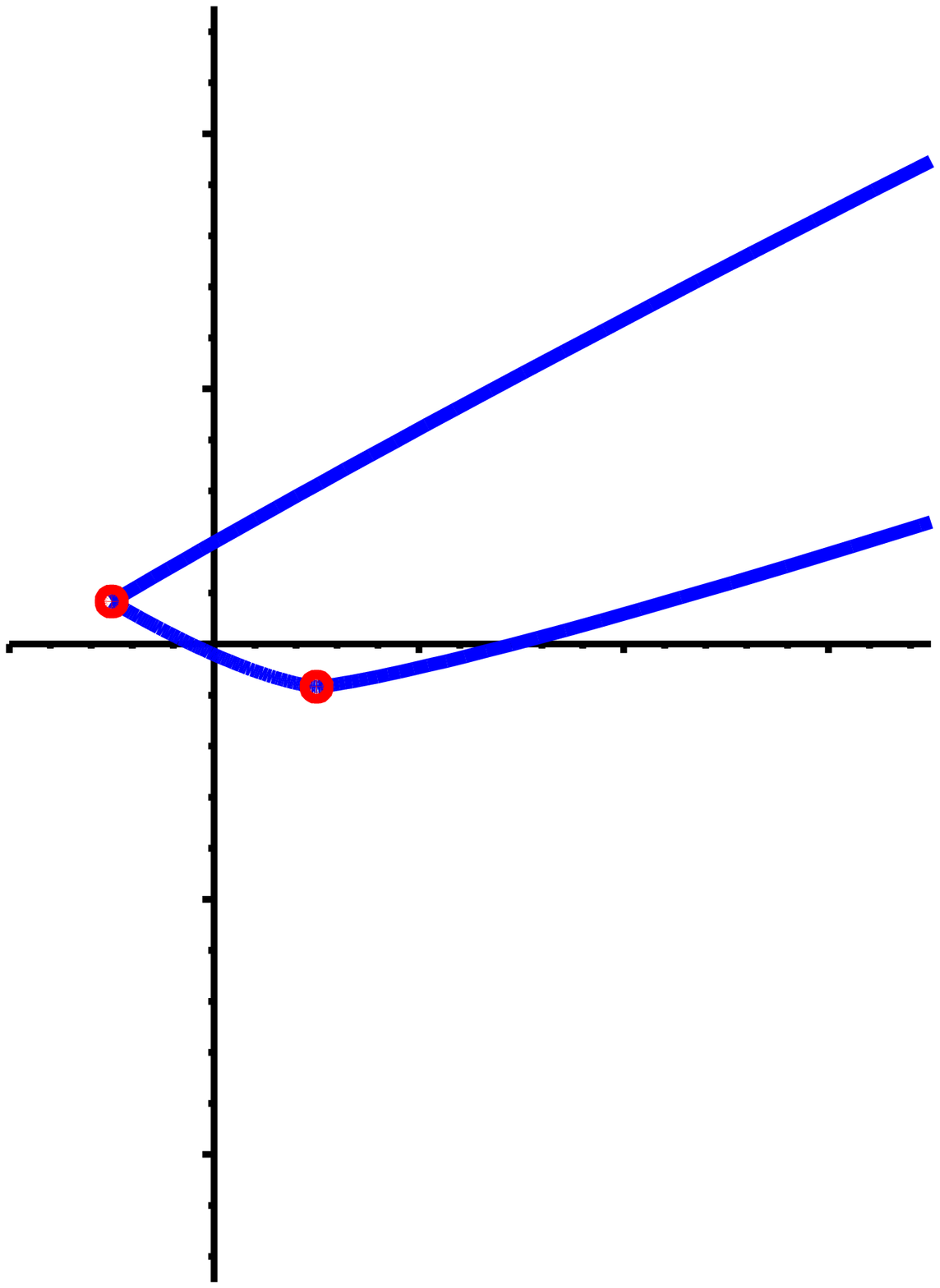}
    \includegraphics[width=2.2cm]{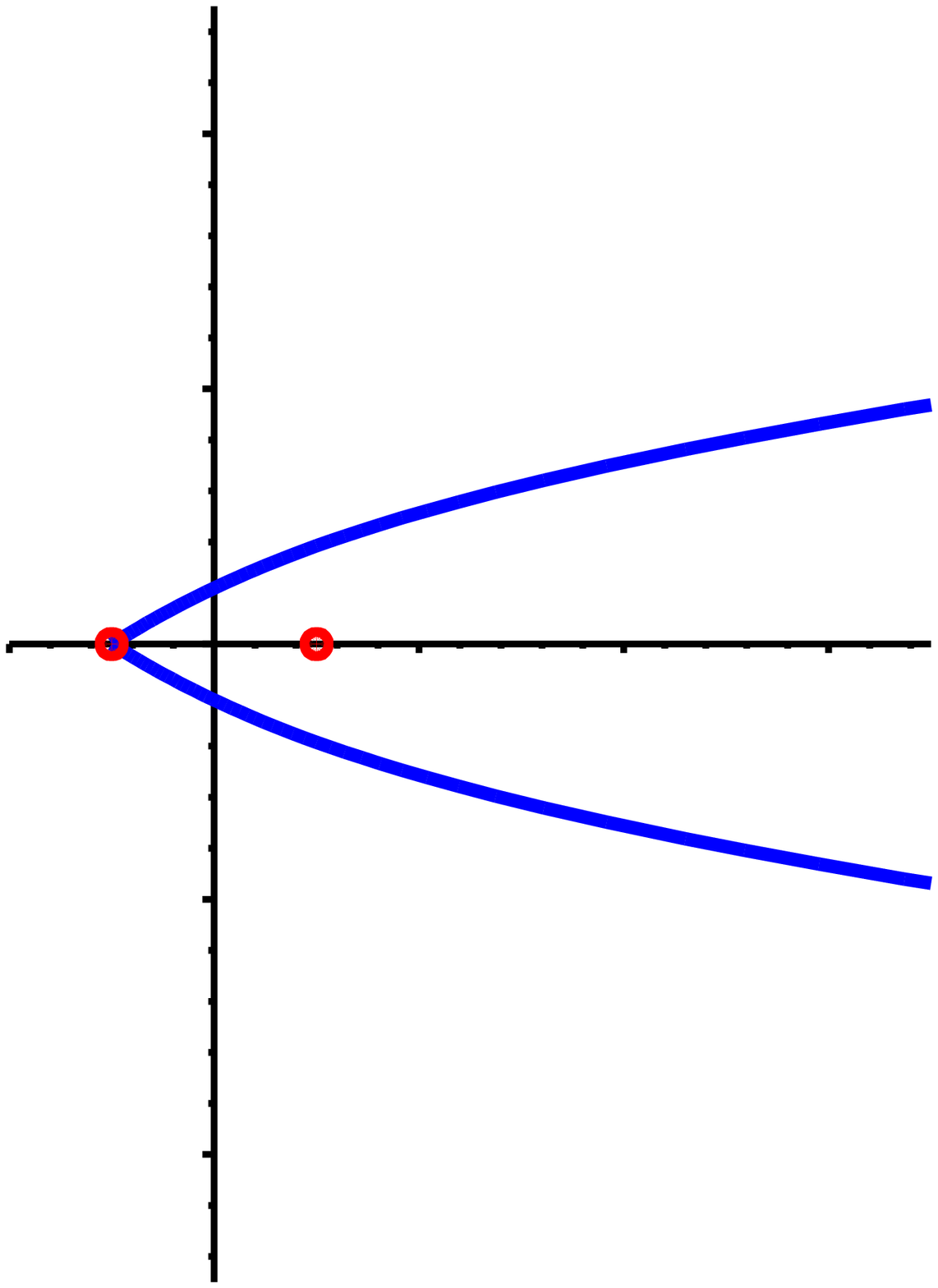}
    \includegraphics[width=2.2cm]{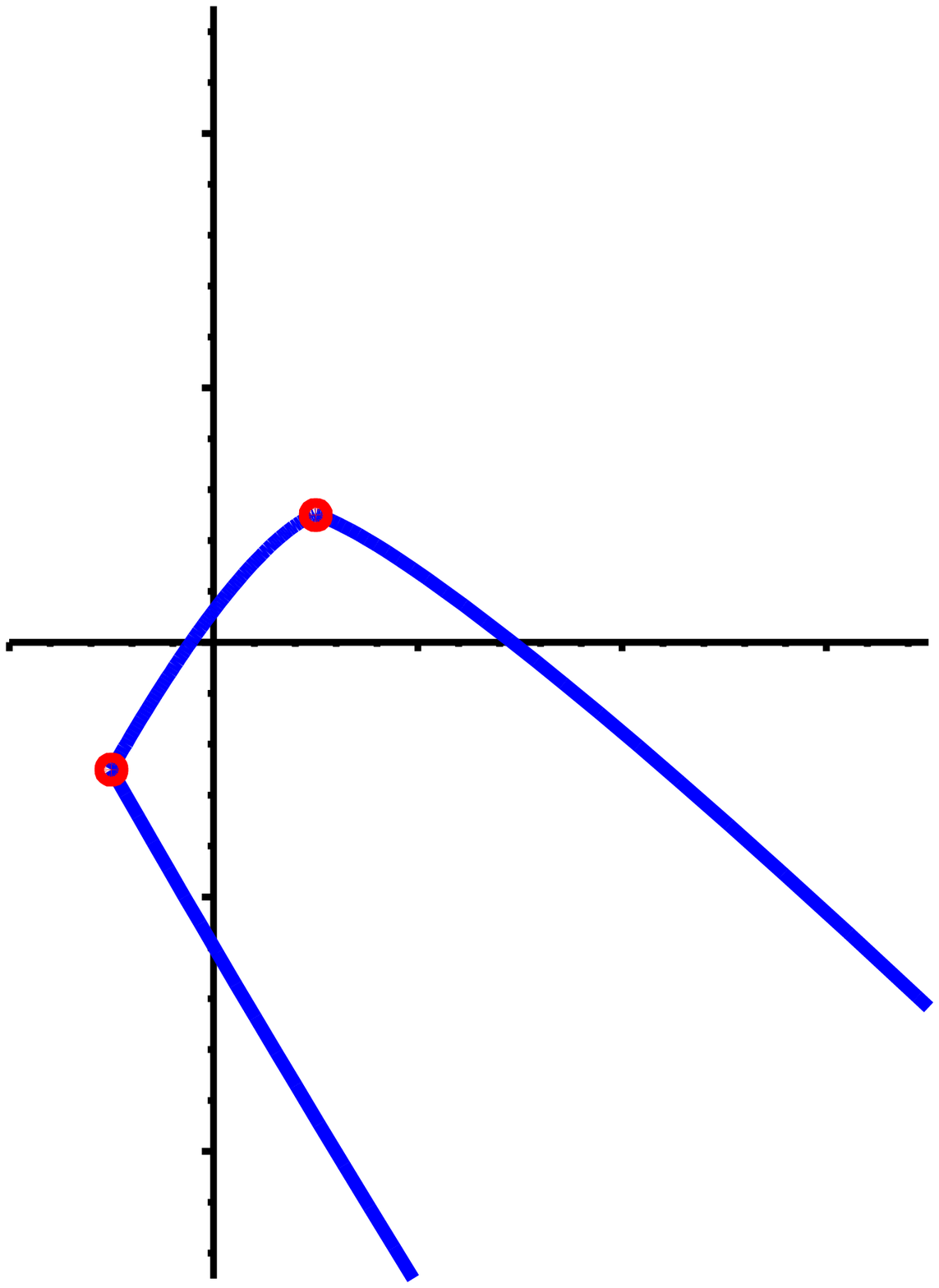}
    \includegraphics[width=2.2cm]{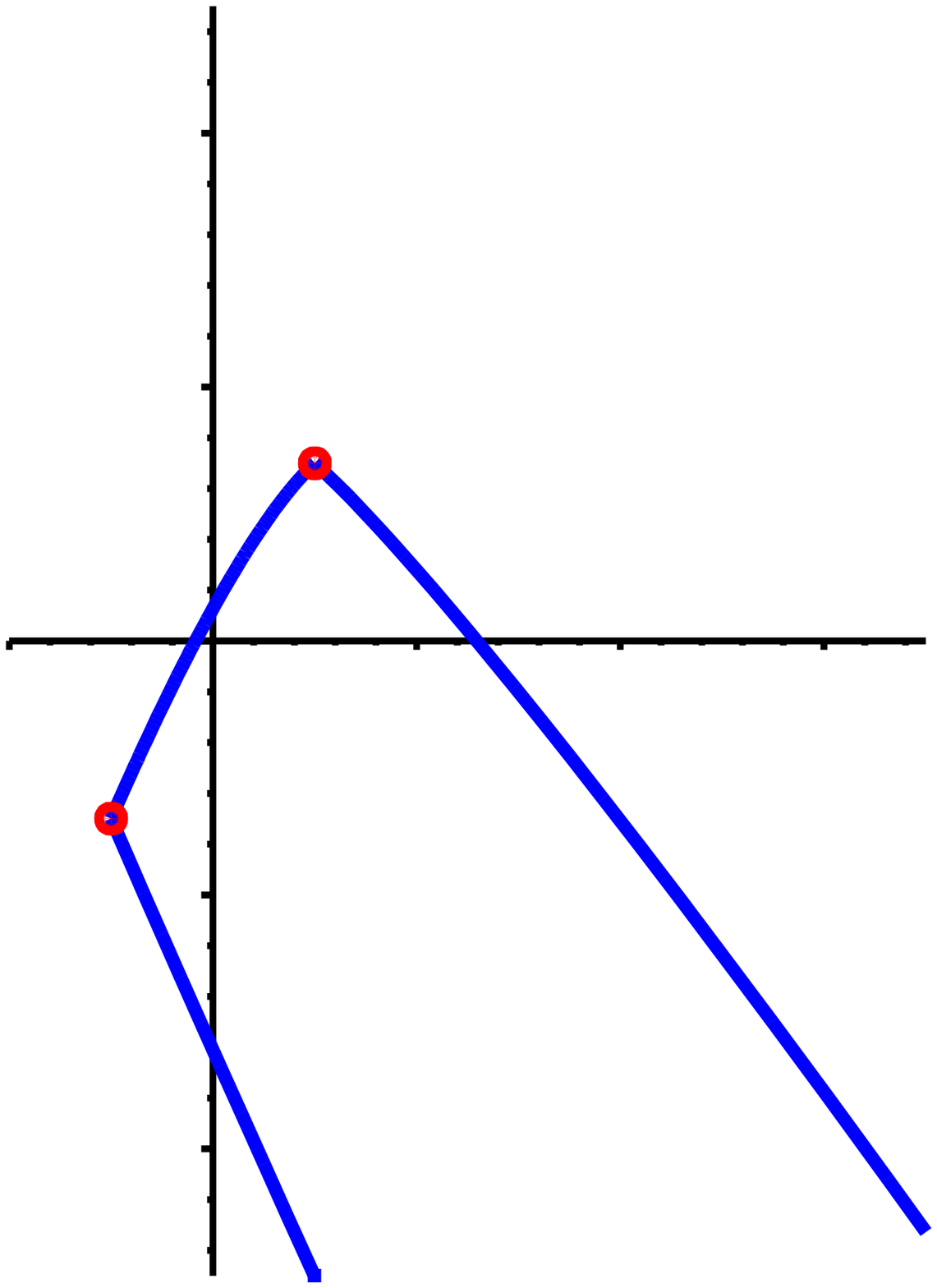}
    
    \emph{Corresponding generalised polytopes:}
    
    \includegraphics[width=2.2cm]{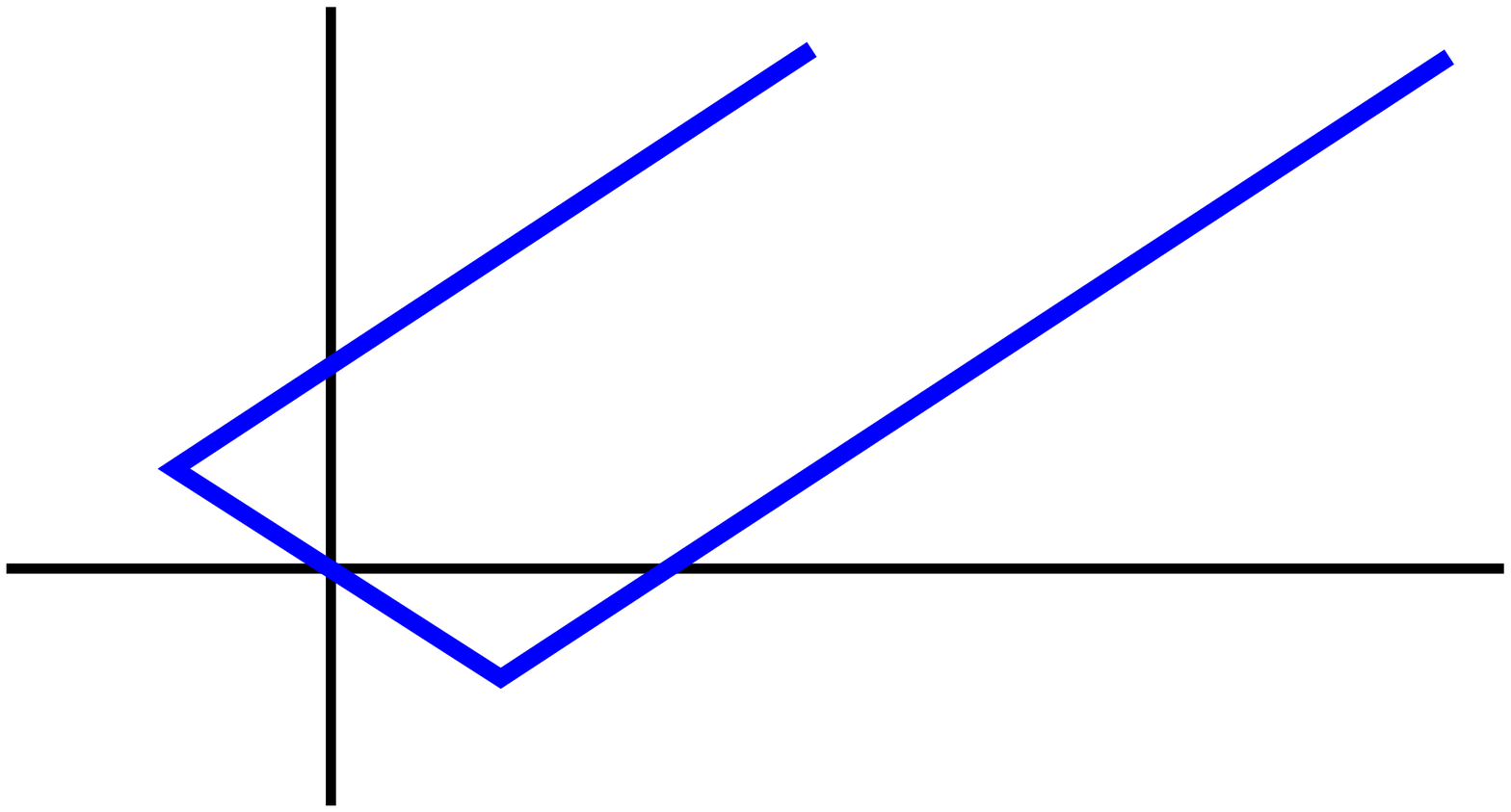}
    \includegraphics[width=2.2cm]{polytope+1}
    \fbox{\begin{minipage}[c]{2.2cm} \vspace{0pt}
        
        \includegraphics[width=\textwidth]{polytope+1}
      
        \centering\emph{and}
      
        \includegraphics[width=\textwidth]{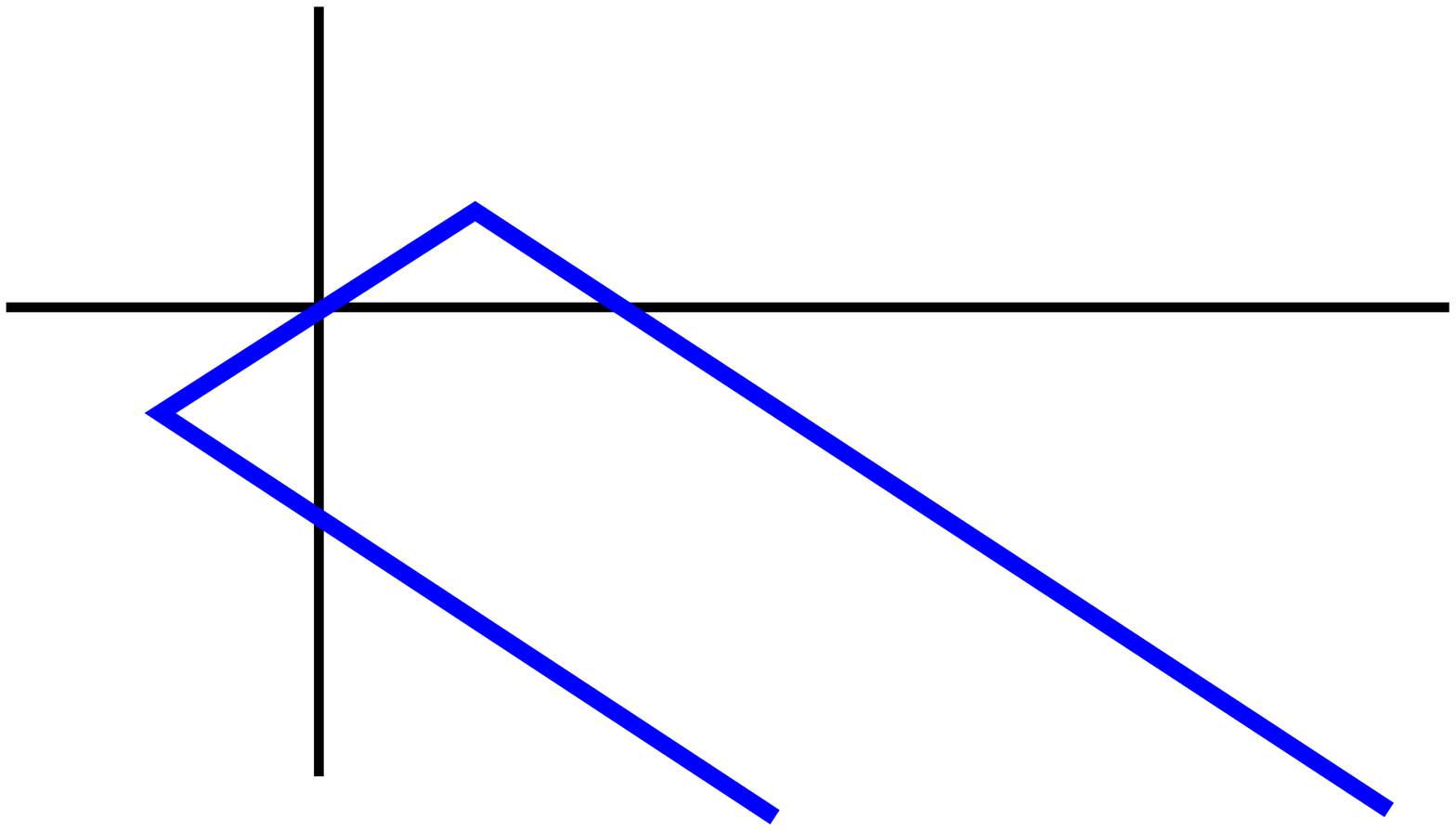}
      \end{minipage}}
    \includegraphics[width=2.2cm]{polytope-1}
    \includegraphics[width=2.2cm]{polytope-1}
    \caption{Bifurcation of the image of the momentum map of example
      \ref{sec:spin}. Here $t=0$, $1/3$, $1/2$, $1$ and $1.2$.}
    \label{fig:moment_t}
  \end{center}
\end{figure}

\begin{demo}
  It is clear that $J$ defines a proper $S¹$ action on $M$. It remains
  to find the singularities of $\Phi$ and compute the spectrum of the
  linearised Hamiltonians; we leave the details to the reader. For
  instance, the two critical points of rank zero are $A_t=(-1,1-2t)$
  and $B_t=(1,-1+2t)$. The spectrum of the linearisation of $H_t$ at
  $A_t$ is composed of two purely imaginary eigenvalues of
  multiplicity two $\pm i\sqrt{5t²-4t+1}$ and hence $A_t$ is always
  elliptic-elliptic, whereas the spectrum at $B_t$ is composed of two
  eigenvalues of multiplicity two $\pm \sqrt{-3t²+4t-1}$, which are
  real if and only if $t\in[1/3,1]$.  In each 2-dimensional eigenspace
  the eigenvalues of $J$ are $\pm i$. Hence $B_t$ is elliptic-elliptic
  for $t<1/3$ and $t>1$ and focus-focus for $t\in(1/3,1)$.
  \begin{figure}[htbp]
    \begin{center}
      \input{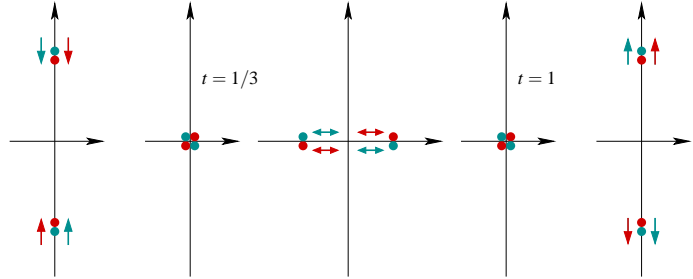}
      \caption{bifurcation of the spectrum of the linearisation of
        $H_t$ at $B_t$}
      \label{fig:bifurc}
    \end{center}
  \end{figure}
\end{demo}

For $t\in(1/3,1)$ we have two generalised polygons. Notice that since
$\Phi_t$ depends continuously on $t$ while the polygons are rational
and hence locally constant, they actually don't depend on
$t\in(1/3,1)$. This of course if also a consequence of the description
in terms of fixed point data. At the south pole (elliptic-elliptic
point) the isotropy weights for $J$ are $(1,1)$ and at the north pole
(focus-focus point) the isotropy weights are $(1,-1)$. We deduce that
the generalised polygons are the one in figure~\ref{fig:moment_0} and
its mirror image with respect to the horizontal axis.

\section{Final remarks}

The construction of the moment polygons for semi-toric momentum maps
was originally motivated by a question of \zhilin{} about the
redistribution of semiclassical eigenvalues in one-parameter families
of quantum Hamiltonian systems. Some hints were given in the very
interesting article~\cite{sadovski-zhilinski}, were the
example~\ref{sec:s2s2} mentioned above was studied from different
viewpoints. In an article in preparation~\cite{san-redistribution} I
give an answer to \zhilin's question in the semi-toric framework.
The moment polygons are a very natural and efficient tool for proving
and stating the result. Roughly speaking, it is shown first using a
global version of Bohr-Sommerfeld rules that the number of eigenvalues
in each ``polyad'' is given in terms of the Duistermaat-Heckman
measure for the Hamiltonian $H$. Secondly, the bifurcation of the
system as the parameter varies is interpreted in terms of an action of
the group $G$ on the initial moment polygon, which gives a geometric
formula for the variation of the Duistermaat-Heckman measure.

Finally I would like to point out that I did not consider in this
article ``inverse questions'' such as which polygons can show up and
to what extent a given class of polygons determines the symplectic
manifold with momentum map $\Phi$. I hope to return on these problems
in a future article, using the invariants of focus-focus foliations of
\cite{san-semi-global}.  However it is easy to see using the
classification by Karshon~\cite{karshon-S1} that in case $M$ is
compact, a given polygon uniquely determines $M$ with the
$S¹$-momentum map $J$. In particular this shows that $M$ always admit
a Kähler structure. But it is not always possible to find a $\T²$
momentum map extending $J$. In view of \zhilin's problem this issue
is not particularly interesting because the initial Lagrangian
foliation would in general be completely different from the toric one
that one could possibly construct (focus-focus leaves do not appear in
toric foliations). For instance it is true that $S²\times S²$ (in
example~\ref{sec:s2s2}) is toric, but this does not help understanding
the redistribution problem, whereas the polygons of
theorem~\ref{theo:polygon} contain all the information we need.

\section*{Acknowledgements}

This article was originally strongly motivated by lively discussions
with Boris \zhilin. A first version was developed under the
auspices of regretted European network MASIE. Later on I had the
chance to get many useful comments by Margaret Symington, which
improved a lot the manuscript. I wish therefore to express all my
gratitude to Boris and Margaret, and to Mark Robert, responsible for
MASIE business.

\bibliographystyle{plain}%
\bibliography{bibli}

\end{document}